\documentclass[generic]{imsart}
\usepackage{amsmath, amssymb,amsfonts,bbm,verbatim,multirow}
\usepackage{epic,eepic,psfrag,epsfig}
\usepackage[sf,bf,SF,footnotesize]{subfigure}
\usepackage{graphicx}
\usepackage{amsthm}
\usepackage{amscd}
\usepackage{epsfig}
\usepackage{natbib}
\usepackage{enumerate}
\usepackage{array}
\usepackage[small]{caption}
\usepackage{url}
\usepackage[left=1.5in,right=1.5in]{geometry}

\newcommand{\FF}{{\mathbb F}}

\newcommand{\NN}{{\mathbb N}}
\newcommand{\RR}{{\mathbb R}}

\newcommand{\UU}{{\mathbb U}}

\DeclareMathOperator*{\argmin}{argmin}

\newtheorem{thm}{Theorem}[section]

\newtheorem{lem}[thm]{Lemma}

\newtheorem{prop}[thm]{Proposition}
\newtheorem{cor}[thm]{Corollary}
\newenvironment{prooftitle}[1]{{\noindent \textsc{Proof #1}}\\}

\numberwithin{equation}{section}

\begin{document}

\begin{frontmatter}

\title{On Convex Least Squares Estimation \\ when the Truth is Linear}
\runtitle{Convex Least Squares Estimation}

\begin{aug}
  \author{\fnms{Yining} \snm{Chen}\ead[label=e1]{y.chen@statslab.cam.ac.uk}}
  \and
  \author{\fnms{Jon A.} \snm{Wellner}\ead[label=e2]{jaw@stat.washington.edu}}

  \address{\\University of Cambridge and University of Washington\\
           \printead{e1,e2}}

  \runauthor{Chen and Wellner}

\end{aug}

\begin{abstract}
We prove that the convex least squares estimator (LSE) attains a $n^{-1/2}$ pointwise rate of convergence
in any region where the truth is linear. In addition, the asymptotic distribution can be characterized
by a modified invelope process. Analogous results hold when one uses the derivative of the convex LSE 
to perform derivative estimation. These asymptotic results facilitate a new consistent testing procedure 
on the linearity against a convex alternative.
Moreover, we show that the convex LSE adapts to the optimal rate 
at the boundary points of the region where the truth is linear, up to a log-log factor. 
These conclusions are valid in the context of both density estimation 
and regression function estimation. 
\end{abstract}



\tableofcontents

\end{frontmatter}
\newpage

\section{Introduction}
Shape-constrained estimation has received much attention recently. The attraction is
the prospect of obtaining automatic nonparametric estimators with no smoothing
parameters to choose. Convexity is among the popular shape constraints that are of both  
mathematical and practical interest. \citet{GJW2001b} show that under the convexity
constraint, the least squares estimator (LSE) can be used to estimate both a density and 
a regression function. For density estimation, 
they showed that the LSE $\hat{f}_n$ of the true convex density $f_0$ converges pointwise at 
a $n^{-2/5}$ rate under certain assumptions.
The corresponding asymptotic distribution can be characterized via a so-called
``invelope'' function investigated by \citet{GJW2001a}. In the regression setting, 
similar results hold for the LSE $\hat{r}_n$ of the true regression function $r_0$. 

However, in the development of their pointwise asymptotic theory, it is required that $f_0$ (or $r_0$) has positive 
second derivative in a neighborhood of the point to be estimated. 
This assumption excludes certain convex functions that may be of practical value. 
Two further scenarios of interest are given below:
\begin{itemize}
\item[1.] At the point $x_0$, the $k$-th derivative $f_0^{(k)}(x_0) = 0$ (or $r_0^{(k)}(x_0) = 0$) for $k = 2,3,\ldots,2s-1$ and $f_0^{(2s)}(x_0) > 0$ (or $r_0^{(2s)}(x_0) > 0$), where $s$ is an integer greater than one;
\item[2.] There exists some region $[a,b]$ on which $f_0$ (or $r_0$) is linear.
\end{itemize}
Scenario~1 can be handled using techniques developed in \citet{BRW2009}. The
aim of this manuscript is to provide theory in the setting of Scenario~2.

We prove that for estimation of a convex density when Scenario~2 holds, 
at any fixed point $x_0 \in (a,b)$, 
the LSE $\hat{f}_n(x_0)$ converges pointwise to 
$f_0(x_0)$ at a $n^{-1/2}$ rate. Its (left or right) derivative $\hat{f}_n'(x_0)$
converges to $f_0'(x_0)$ at the same rate. 
The corresponding asymptotic distributions are characterized using 
a modified invelope process. More generally, we show that the processes $\big\{\sqrt{n}(\hat{f}_n(x) - f_0(x)) : x \in [a+\delta,b-\delta]\big\}$ 
and $\big\{\sqrt{n}(\hat{f}_n'(x) - f_0'(x)) : x \in [a+\delta,b-\delta]\big\}$ both converge weakly for any $\delta > 0$.
We remark that unlike the case of \citet{GJW2001b}, these processes are not local, because intuitively speaking, the linear behavior of the true function is non-local, which requires linearity on an interval with positive length. In addition, there does not exist a 
``universal'' distribution of $\hat{f}_n$ on $(a,b)$, 
i.e. the pointwise limit distributions at different points are in general different. 
In addition, we show that the derived asymptotic processes can be used in a more practical setting where one would like to 
perform tests on the linearity against a convex alternative.
Moreover, we study the adaptation of the LSE $\hat{f}_n$ at the boundary points of the linear region (e.g. $a$ and $b$). 
Note that the difficulty level of estimating $f_0(a)$ and $f_0(b)$ depends on the behavior
of $f_0$ outside $[a,b]$. Nevertheless, we show that $\hat{f}_n(a)$ (or $\hat{f}_n(b)$) converges to $f_0(a)$ (or $f_0(b)$)
at the minimax optimal rate up to a negligible factor of $\sqrt{\log \log n}$.
Last but not least, we establish the analogous rate and asymptotic distribution results
for the LSE $\hat{r}_n$ in the regression setting.
 
Our study yields a better understanding of the adaptation 
of the LSE in terms of pointwise convergence under the convexity constraint. 
It is also one of the first attempts to quantify 
the behavior of the convex LSE at non-smooth points. 
When the truth is linear, the minimax optimal $n^{-1/2}$ pointwise rate is indeed achieved by the LSE
on $(a,b)$. The optimal rate at the boundary points $a$ and $b$ is also achievable by the LSE up to a
log-log factor. These results can also be applied to the case where the true function consists of multiple linear components.
Furthermore, our results can be viewed as an intermediate stage for the development 
of theory under misspecification. Note that linearity is regarded as the boundary case
of convexity: if a function is non-convex, then its projection to the class of convex 
functions will have linear components. We conjecture that the LSE in these 
misspecified regions converges at an $n^{-1/2}$ rate, with the asymptotic distribution
characterized by a more restricted version of the invelope process. 
More broadly, we expect that this type of behavior will be seen in situations of other shape 
restrictions, such as log-concavity for $d = 1$ \citep{BRW2009} and $k$-monotonicity 
\citep{BalabdaouiWellner2007}.

The LSE of a convex density function was first studied by \citet{GJW2001b}, where its consistency
and some asymptotic distribution theory were provided. On the other hand, the idea of using the LSE for 
convex regression function estimation dates back to \citet{Hildreth1954}. 
Its consistency was proved by \citet{HansonPledger1976}, with some rate results given 
in \citet{Mammen1991}. In this manuscript, for the sake of mathematical convenience, 
we shall focus on the non-discrete version discussed by \citet{BalabdaouiRufibach2008}.
See \citet{GJW2008} for the computational aspects of all the above-mentioned LSEs.

There are studies similar to ours regarding other shape restrictions. See 
Remark~2.2 of \citet{Groeneboom1985} and
\citet{CarolanDykstra1999} in the context of decreasing density function 
estimation (a.k.a. Grenander estimator) when the truth is flat and \citet{Balabdaoui2014} with regard to 
discrete log-concave distribution estimation when the true distribution is geometric. 
In addition, \citet{GroeneboomPyke1983} studied the Grenander estimator's global behavior in the $L^2$ norm under
the uniform distribution and gave a connection to statistical problems involving combination of $p$-values and two-sample rank statistics, 
while \citet{Cator2011} studied the adaptivity of the LSE in monotone regression (or density) function estimation. 
For estimation under misspecification of various shape constraints, we point readers to 
\citet{Jankowski2014}, \citet{CuleSamworth2010}, \citet{DSS2011}, 
\citet{ChenSamworth2013}, and \citet{BJRP2013}. More recent developments on global
rates of the shape-constrained methods can be found in 
\citet{GuntuboyinaSen2015}, \citet{DossWellner2016}, and \citet{KimSamworth2014}.
See also \citet{Meyer2013} and \citet{ChenSamworth2014} where an additive structure is imposed
in shape-constrained estimation in the multidimensional setting. 

The rest of the paper is organized as follows: in Section~\ref{Sec:density}, we study the 
behavior of the LSE for density estimation. In particular, for notational convenience, we first 
focus on a special case where the true density function $f_0$ is taken to be triangular.
The convergence rate and asymptotic distribution are given in Section~\ref{Sec:densityrate} and Section~\ref{Sec:densityasym}. 
Section~\ref{Sec:test} demonstrates the practical use of our asymptotic results, 
where a new consistent testing procedure on the linearity against a convex alternative is proposed. 
More general cases are handled later in Section~\ref{Sec:densityext} based on the ideas illustrated 
in Section~\ref{Sec:densityspecial}. Section~\ref{Sec:densityadapt}
discusses the adaptation of the LSE at the boundary points. 
Analogous results with regard to regression function estimation are presented in Section~\ref{Sec:regress}. 
Some proofs, mainly on the existence and uniqueness of a limit process and the adaptation of the LSE,  
are deferred to the appendices.

\section{Estimation of a density function}
\label{Sec:density}
Given $n$ independent and identically distributed (IID) observations from a density function $f_0: [0,\infty) \rightarrow [0,\infty)$. 
Let $F_0$ be the corresponding distribution function (DF).
In this section, we denote the convex cone of all non-negative continuous convex and integrable functions on $[0,\infty)$ by $\mathcal{K}$.
The LSE of $f_0$ is given by 
\[
\hat{f}_n  = \argmin_{g \in \mathcal{K}} \Big(\frac{1}{2}\int_0^{\infty} g(t)^2dt - \int_0^{\infty} g(t) d \mathbb{F}_n(t) \Big),
\]
where $\mathbb{F}_n$ is the empirical distribution function of the observations. Furthermore, we denote the DF of $\hat{f}_n$ by $\hat{F}_n$. 

Throughout the manuscript, without specifying otherwise, the derivative of a convex function can be interpreted as either its left derivative or its right derivative.

\subsection{A special case}
\label{Sec:densityspecial}

To motivate the discussion, as well as for notational convenience, we take $f_0 (t) = 2(1-t) \mathbf{1}_{[0,1]}(t)$ in Section~\ref{Sec:densityrate} and Section~\ref{Sec:densityasym}. 
Extensions to other convex density functions are presented in Section~\ref{Sec:densityext}.

Consistency of $\hat{f}_n$ over $(0,\infty)$ in this setting can be found in \citet{GJW2001b}. \citet{Balabdaoui2007} has shown how $\hat{f}_n$ can be used to provide consistent estimators of $f_0(0)$ and $f_0'(0)$. 

Here we concentrate on the LSE's rate of convergence and asymptotic distribution. 

\subsubsection{Rate of convergence}
\label{Sec:densityrate}
The following proposition shows that the convergence rate at any interior point where $f_0$ is linear is $n^{-1/2}$, thus achieving the optimal rate given in Example~1 of \citet{CaiLow2015}.
\begin{thm}[Pointwise rate of convergence]
\label{Thm:LS_triangular_rate}
Suppose $f_0 (t) = 2(1-t) \mathbf{1}_{[0,1]}(t)$. Then for any fixed $x_0 \in (0,1)$, 
\[
	|\hat{f}_n(x_0) - f_0(x_0)| = O_p(n^{-1/2}).
\]
\end{thm}

\begin{prooftitle}{of Theorem~\ref{Thm:LS_triangular_rate}}
A key ingredient of this proof is the version of Marshall's lemma in this setting \citep[Theorem 1]{DRW2007}, which states that 
\begin{align}
\label{Eq:marshall}
	\lVert \hat{F}_n - F_0 \rVert_{\infty} \le 2 \lVert \mathbb{F}_n - F_0 \rVert_{\infty},
\end{align}
where $\lVert \cdot \rVert_\infty$ is the uniform norm.

Let $c = \hat{f}_n(x_0) - f_0(x_0)$. Two cases are considered here: (a) $c > 0$ and (b) $c < 0$. 

In the first case, because $\hat{f}_n$ is convex, one can find a supporting hyperplane of $\hat{f}_n$ passing through $(x_0,f_0(x_0)+c)$ such that
\[
	\hat{f}_n(y) \ge  k (y - x_0) + f_0(x_0) + c, \quad \mbox{ for any } y \in[0,\infty),
\]
where $k$ is a negative slope. If $k \le -2$, then $\hat{F}_n(x_0) - F_0(x_0) \ge cx_0 > 0$. Otherwise, $\hat{F}_n(x_0) - F_0(x_0) < - c(1-x_0) < 0$. Figure~\ref{fig:greater} explains the above inequalities graphically. Consequently, $|\hat{F}_n(x_0) - F_0(x_0)| \ge c \, \min(x_0,1-x_0)$.

\begin{figure}[ht!]
  \centering
  $\begin{array}{c c} 
  \includegraphics[scale=0.37]{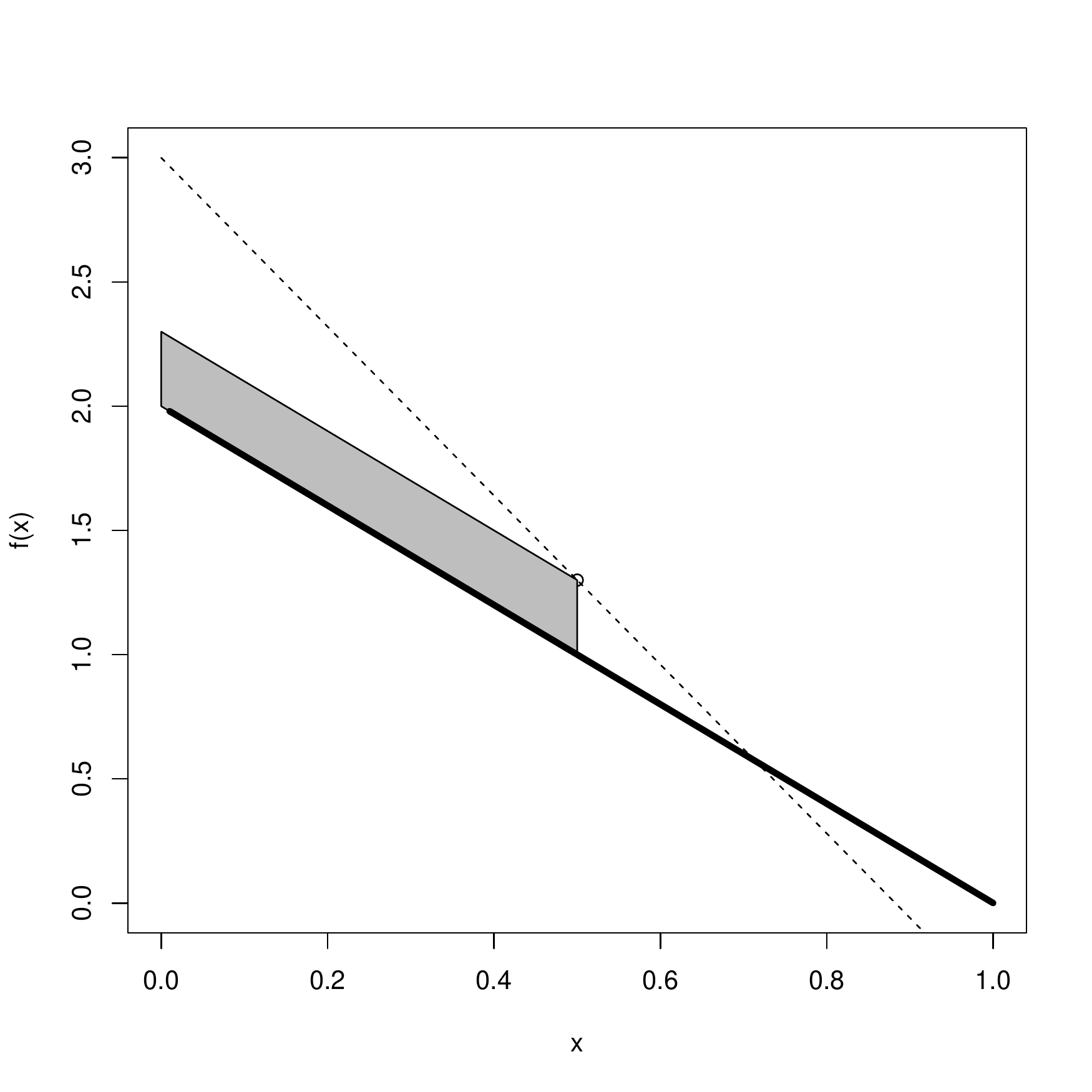} &
  \includegraphics[scale=0.37]{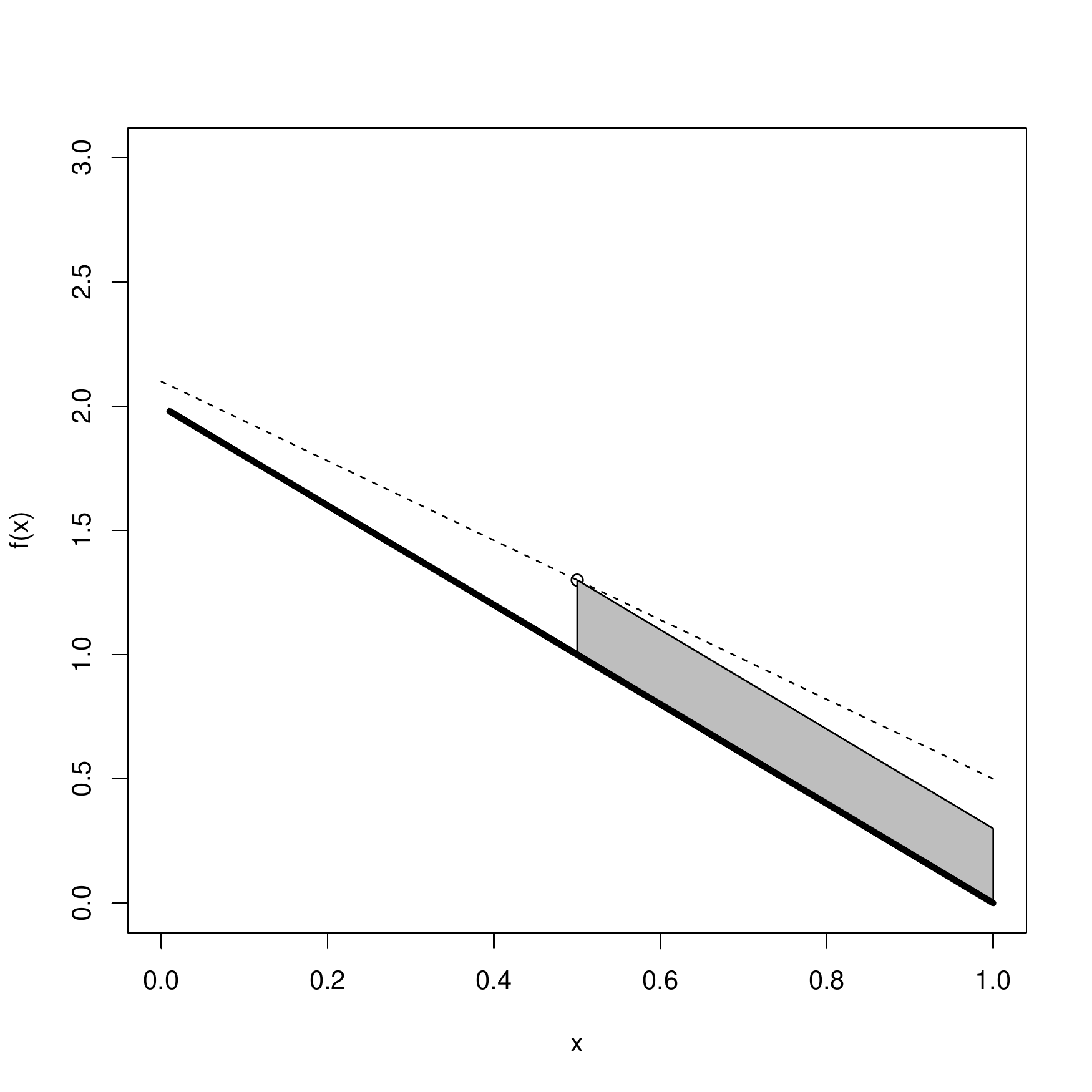} \\
  \mathrm{(A)} & \mathrm{(B)} 
  \end{array}$
\caption{The case of $c > 0$: (A) illustrates the situation of $k \le -2$, where $\hat{F}_n(x_0) - F_0(x_0) \ge cx_0$; (B) illustrates the situation of $k > -2$, where $\hat{F}_n(x_0) - F_0(x_0) < - c(1-x_0)$. In both scenarios, the dashed curve represents the supporting hyperplane of the LSE, while the thick solid line represents the true density $f_0$.}
\label{fig:greater}
\end{figure}
\begin{figure}[ht!]
  \centering
  $\begin{array}{c c} 
  \includegraphics[scale=0.37]{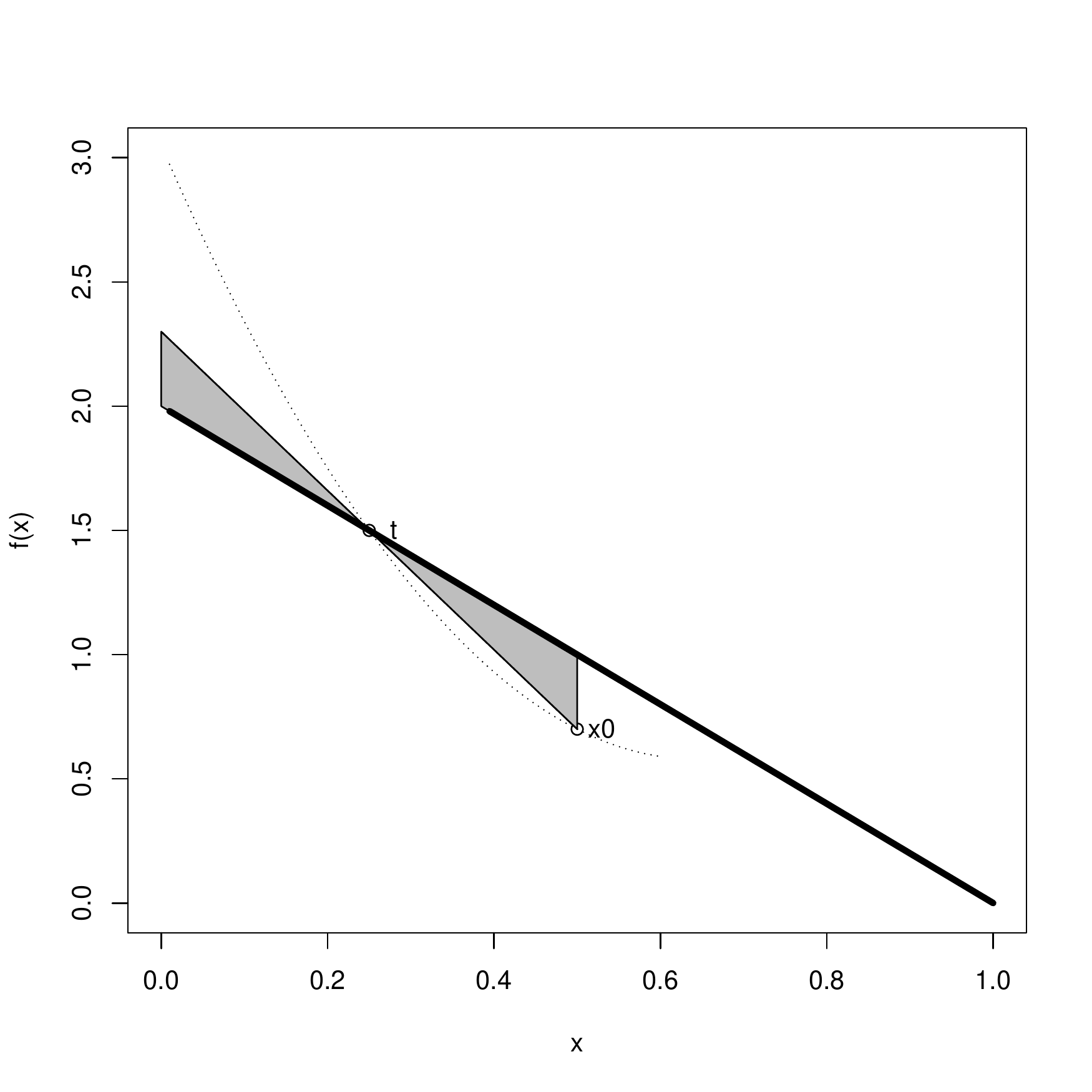} &
  \includegraphics[scale=0.37]{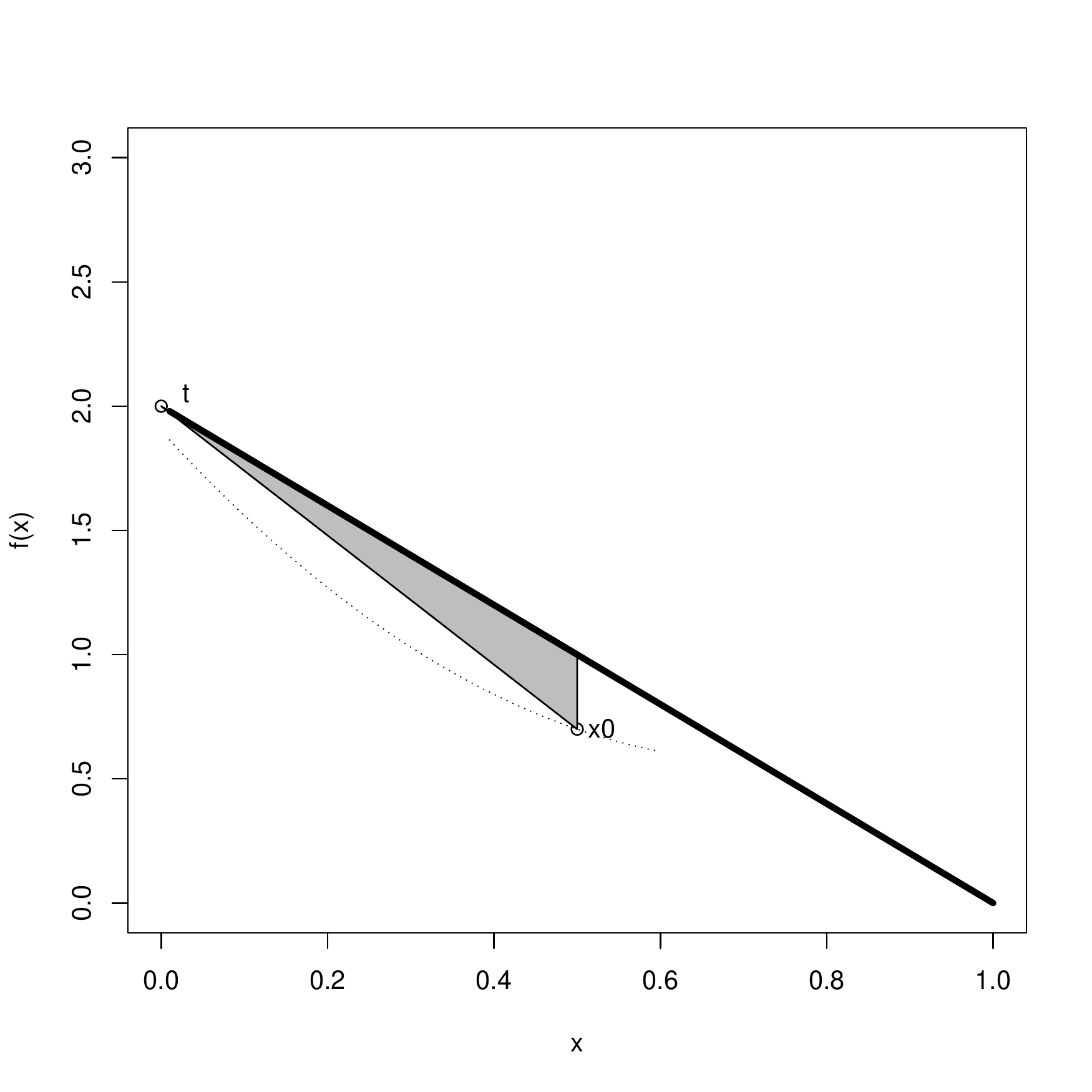} \\
  \mathrm{(A)} & \mathrm{(B)} 
  \end{array}$
 \caption{The case of $c < 0$: (A) illustrates the inequalities when $\hat{f}_n$ and $f_0$ has an intersection on $[0,x_0)$; the situation where $\hat{f}_n$ and $f_0$ has no intersection on $[0,x_0)$ is represented in (B). In both scenarios, the dotted curve represents the LSE and the thick solid line represents the true density $f_0$.}
\label{fig:smaller}
\end{figure}

In the second case, if we can find $0 \le t < x_0$ such that $\hat{f}_n$ and $f_0$ intersect at the point $(t,2-2t)$, then 
\begin{align*}
	\hat{F}_n(t) - F_0(t) &\ge \frac{1}{2}t\frac{|c|t}{x_0-t} \ge 0; \\
	\big(\hat{F}_n(x_0) - \hat{F}_n(t)\big) - \big(F_0(x_0) - F_0(t)\big) &\le -\frac{1}{2}|c|(x_0-t) < 0.
\end{align*}
Otherwise, by taking $t=0$, we can still verify that the above two inequalities hold true. Figure~\ref{fig:smaller} illustrates these findings graphically. 

Therefore, 
\begin{align*}
	2\lVert \hat{F}_n - F_0 \rVert_\infty &\ge \max\Big\{|\hat{F}_n(t) - F_0(t)|,\big|\big(\hat{F}_n(x_0) - F_0(x_0)\big) - \big(\hat{F}_n(t) - F_0(t)\big)\big|\Big\} \\
	&\ge \frac{|c|}{2}\max\Big(\frac{t^2}{x_0-t},x_0-t\Big) \ge \frac{|c|x_0}{4},
\end{align*}
where the last inequality uses the fact that for any $t \in [0,x_0)$, $t^2/(x_0-t)$ is an increasing function of $t$ while $x_0-t$ is a decreasing function of $t$.

By (\ref{Eq:marshall}), we have that
\begin{align*}
	O_p(n^{-1/2}) = 4 \lVert \mathbb{F}_n - F_0 \rVert_{\infty} \ge 2 \lVert \hat{F}_n - F_0 \rVert_{\infty} \ge \min\left( 2x_0,2-2x_0,\frac{x_0}{4}\right) |c|.
\end{align*}

It then follows that $c = O_p(n^{-1/2})$, as desired. \hfill $\Box$
\end{prooftitle}

\begin{cor}[Uniform rate of convergence]
\label{Cor:LS_triangular_uniform_rate}
For any $0 < \delta \le 1/2$, 
\[
	\sup_{x \in [\delta, 1-\delta]} |\hat{f}_n(x) - f_0(x)| = O_p(n^{-1/2}).
\]
\end{cor}

\begin{prooftitle}{of Corollary~\ref{Cor:LS_triangular_uniform_rate}}
A closer look at the proof of Theorem~\ref{Thm:LS_triangular_rate} reveals that we have
\begin{align*}
	O_p(n^{-1/2}) = 4 \lVert \mathbb{F}_n - F_0 \rVert_{\infty} &\ge \min\left( 2x,2-2x,\frac{x}{4}\right) |\hat{f}_n(x) - f_0(x)|
\end{align*}
simultaneously for every $x \in [\delta, 1-\delta]$. Therefore, 
\[
	O_p(n^{-1/2}) = \inf_{x \in [\delta, 1-\delta]}\min\left( 2x,2-2x,\frac{x}{4}\right) \sup_{x \in [\delta, 1-\delta]}|\hat{f}_n(x) - f_0(x)|.
\]
Consequently, $\sup_{x \in [\delta, 1-\delta]}|\hat{f}_n(x) - f_0(x)|$ is $O_p(n^{-1/2})$.
\hfill $\Box$
\end{prooftitle}

Let $\hat{f}^-_n$ and $\hat{f}^+_n$ denote respectively the left and right derivatives of $\hat{f}_n$. The same convergence rate also applies to these derivative estimators.
\begin{cor}[Uniform rate of convergence: derivatives]
\label{Cor:LS_triangular_uniform_rate_derivative}
For any $0 < \delta \le 1/2$, 
\[
	\sup_{x \in [\delta, 1-\delta]} \max \left(|\hat{f}^-_n(x) - f_0'(x)|, |\hat{f}^+_n(x) - f_0'(x)| \right) = O_p(n^{-1/2}).
\]
\end{cor}

\begin{prooftitle}{of Corollary~\ref{Cor:LS_triangular_uniform_rate_derivative}}
By the convexity of $\hat{f}_n$,
\begin{align*}
	&\sup_{x \in [\delta, 1-\delta]} \max \left(|\hat{f}^-_n(x) - f_0'(x)|, |\hat{f}^+_n(x) - f_0'(x)|\right) \\
	&\le \max \left(|\hat{f}^-_n(\delta)- f_0'(\delta)|, |\hat{f}^+(1-\delta) - f_0'(1-\delta)|\right) \\
	&\le \frac{2}{\delta} \max \Big(\big|\hat{f}_n(\delta) - \hat{f}_n(\delta/2) -  f_0(\delta) + f_0(\delta/2)\big|, \\
	& \qquad \qquad \qquad\big|\hat{f}_n(1- \delta/2) - \hat{f}_n(1-\delta) -  f_0(1- \delta/2) + f_0(1-\delta)\big| \Big) \\
	& = O_p(n^{-1/2}),
\end{align*}
where the final equation follows from Corollary~\ref{Cor:LS_triangular_uniform_rate}.
\hfill $\Box$
\end{prooftitle}

\subsubsection{Asymptotic distribution}
\label{Sec:densityasym}
To study the asymptotic distribution of $\hat{f}_n$, we start by characterizing the limit distribution.
\begin{thm}[Characterization of the limit process]
\label{Thm:triangular_limit}
Let $X(t) = \UU(F_0(t))$ and $Y(t) = \int_0^t X(s) ds$ for any $t \ge 0$, where $\UU$ is a standard Brownian bridge process on $[0,1]$.
Then almost surely (a.s.), there exists a uniquely defined random continuously differentiable function $H$ on $[0,1]$ satisfying the following conditions:
\begin{enumerate}[(1)]
\item $H(t) \ge Y(t)$ for every $t \in [0,1]$;
\item $H$ has convex second derivative on $(0,1)$;
\item $H(0) = Y(0)$ and $H'(0) = X(0)$; 
\item $H(1) = Y(1)$ and $H'(1) = X(1)$; 
\item $\int_0^1 \big(H(t)-Y(t)\big)\, dH^{(3)}(t) = 0$.
\end{enumerate}
The above claim also holds if we instead consider a standard Brownian motion $\tilde{X}(t)$ on $[0,1]$ and replace $(X,Y,H)^T$ by $(\tilde{X},\tilde{Y},\tilde{H})^T$, where $\tilde{Y}(t) = \int_0^t \tilde{X}(s) ds$, and where $\tilde{H}$ is different from $H$. 
\end{thm}

\begin{figure}[ht!]
  \centering
  \includegraphics[height=8cm, width=14cm, trim=0 220 0 190]{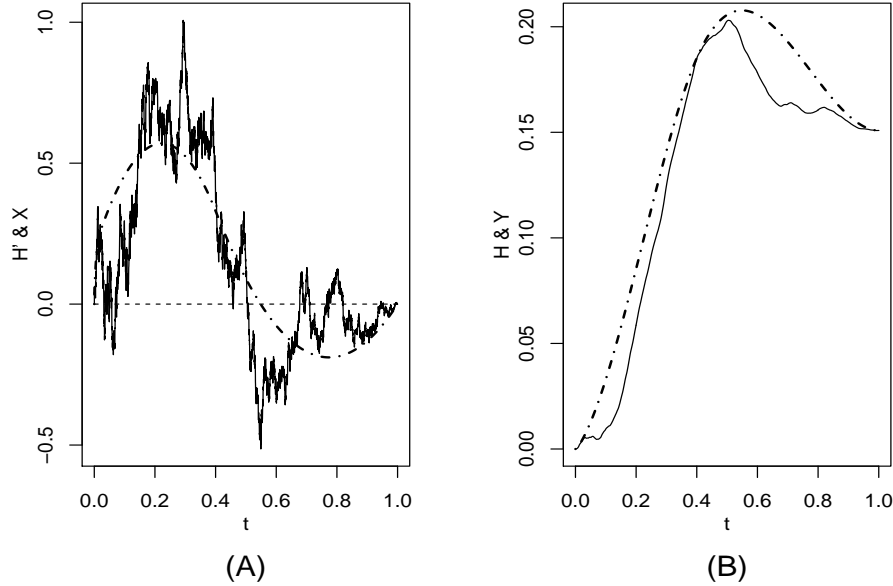} 
\caption{A typical realization of $X(t)$, $Y(t)$, $H'(t)$ and $H(t)$ in Theorem~\ref{Thm:triangular_limit} with $X(t) = \UU(F_0(t))$ and $F_0$ the distribution function with triangular density $f_0(t) = 2(1-t) \mathbf{1}_{[0,1]}(t)$: in (A), $X(t)$ and $H'(t)$ are plotted in solid and dash-dotted curves respectively; in (B), the corresponding $Y(t)$ is plotted in solid curve, while the invelope process $H(t)$ is illustrated in dash-dotted curve. }
\label{Fig:XYH}
\end{figure}

Figure~\ref{Fig:XYH} shows a typical realization of $X(t)$, $Y(t)$, $H'(t)$ and $H(t)$. A detailed construction of the above limit invelope process can be found in Appendix I. Note that our process is defined on a compact interval, so is technically different from the process presented in \citet{GJW2001a} (which is defined on the whole real line). As a result, extra conditions regarding the behavior of $H$ (and $H'$) at the boundary points are imposed here to ensure its uniqueness. 

Other characterization of the limit process is also possible. A slight variant is given below, which proof can also be found in Appendix I.

\begin{cor}[Different characterization]
\label{Cor:triangular_limit}
Conditions (3) and (4) in the statement of Theorem~\ref{Thm:triangular_limit} can be replaced respectively by
\begin{enumerate}
\item[(3')] $\lim_{t \rightarrow 0^+} H^{(2)}(t) = \infty$, and
\item[(4')] $\lim_{t \rightarrow 1^-} H^{(2)}(t) = \infty$.
\end{enumerate}
\end{cor}

Now we are in the position to state our main result of this section.
\begin{thm}[Asymptotic distribution]
\label{Thm:LS_triangular_asymptotic}
Suppose $f_0(t) = 2(1-t) \mathbf{1}_{[0,1]}(t)$. Then for any $\delta > 0$, the process
\begin{align*}
\sqrt{n}
\left(\begin{array}{c}
\hat{f}_n (x) - f_0 (x)\\
\hat{f}_n' (x) - f_0' (x) 
\end{array} \right)
\Rightarrow
\left(\begin{array}{c}
H^{(2)}(x) \\
H^{(3)} (x)
\end{array} \right) \mbox{ in }\, \mathcal{C}[\delta, 1-\delta] \times \mathcal{D}[\delta, 1-\delta],
\end{align*}
where $\mathcal{C}$ is the space of continuous functions equipped with the uniform norm, $\mathcal{D}$ is the Skorokhod space, and $H$ is the invelope process defined in Theorem~\ref{Thm:triangular_limit}.

In particular, for any $x_0 \in (0,1)$,
\begin{align}
\label{Eq:LS_pointwise_converg}
\sqrt{n} \Big(\hat{f}_n(x_0) - f_0(x_0), \hat{f}_n'(x_0) - f_0'(x_0)\Big)^T \stackrel{d}{\rightarrow}  \Big(H^{(2)}(x_0),H^{(3)}(x_0)\Big)^T.
\end{align}
\end{thm}
\begin{prooftitle}{of Theorem~\ref{Thm:LS_triangular_asymptotic}}
Before proceeding to the proof, we first define the following processes on $[0, \infty)$: 
\begin{align}
\notag	X_n(t) &= n^{1/2}  \Big(\FF_n(t) - F_0(t)\Big); \\
\notag	Y_n(t) &= n^{1/2} \int_{0}^t \Big(\FF_n(s) - F_0(s)\Big) ds ;\\
\label{Eq:Def_H_n}	\hat{H}_n(t) &= n^{1/2} \int_{0}^t \bigg\{\int_{0}^s \Big(\hat{f}_n(v) - f_0(v)\Big) dv \bigg\} ds.
\end{align}
Furthermore, define the set of ``knots'' of a convex function $f$ on $(0,1)$ as 
\begin{align}
\label{Eq:defknots}
	\mathcal{S}(f) = \left\{ t \in (0,1): f(t) < \frac{1}{2} \Big(f(t - \delta) + f(t + \delta)\Big),\; \forall \delta \in \Big(0,\min(1-t,t)\Big) \right\}.
\end{align}
We remark that the above definition of knots can be easily extended to convex functions with a different domain.

By Lemma~2.2 of \citet{GJW2001b}, $\hat{H}_n(t) \ge Y_n(t)$ for $t \in [0,1]$, with equality if $t \in \mathcal{S}(\hat{f}_n)$. In addition, $\int_{0}^1 (\hat{H}_n(t) - Y_n(t))d\hat{H}_n^{(3)}(t) = 0$. 

Now define the space $E_m$ of vector-valued functions ($m \ge 3$) as 
\begin{align*}
	E_m = \mathcal{C}[0,1] \times \mathcal{C}[0,1] \times \mathcal{C}[1/m, 1-1/m] \times \mathcal{D}[1/m, 1-1/m] \times \mathcal{C}[0,1] \times \mathcal{D}[0,1].
\end{align*}
and endow $E_m$ with the product topology induced by the uniform norm on $\mathcal{C}$ (i.e. the space of continuous functions) and Skorokhod metric on $\mathcal{D}$ (i.e. the Skorokhod space). Let $E_m$ be supported by the stochastic process
\begin{align*}
	Z_n = \left(\hat{H}_n, \hat{H}_n', \hat{H}_n^{(2)}, \hat{H}_n^{(3)}, Y_n, X_n \right)^T.
\end{align*}

Corollary~\ref{Cor:LS_triangular_uniform_rate_derivative} entails the tightness of $\hat{H}_n^{(3)}$ in $\mathcal{D}[1/m, 1-1/m]$, while Corollary~\ref{Cor:LS_triangular_uniform_rate} entails the tightness of $\hat{H}_n^{(2)}$ in $\mathcal{C}[1/m, 1-1/m]$. In addition, both $\hat{H}_n$ and $\hat{H}_n'$ are tight in $\mathcal{C}[0,1]$ (via an easy application of Marshall's lemma). Since $X_n$ converges to a Brownian bridge, it is tight in $\mathcal{D}[0,1]$ as well. Finally, the same is also true for $Y_n$. 

Note that $E_m$ is separable. For any subsequence of $Z_n$, by Prohorov's theorem, we can construct a further subsequence $Z_{n_j}$ such that $\{Z_{n_j}\}_j$ converges weakly in $E_m$ to some 
\[
Z_0 = \Big(H_0, H_0', H_0^{(2)}, H_0^{(3)}, Y, X\Big)^T.
\]
In addition, it follows from a diagonalization argument that a further extraction can make $\{Z_{n_j}\}_j$ converge weakly to $Z_0$ in every $E_m$. Using Skorokhod's representation theorem, we can assume without loss of generality that for almost every $\omega$ in the sample space, $Z_{n_j}(\omega) \rightarrow Z_0(\omega)$. Moreover, since $X$ has a continuous path a.s., the convergence of $X_{n_j}$ can be strengthened to $\lVert X_{n_j}(\omega) - X(\omega) \rVert_{\infty} \rightarrow 0$, where $\lVert \cdot \rVert_{\infty}$ is the uniform norm. In the rest of the proof, $\omega$ is suppressed for notational convenience, so depending on the context, $Z_{n_j}$ (or $Z_0$) can either mean a random variable or a particular realization.  

In the following, we shall verify that $H_0$ satisfies all the conditions listed in the statement of Theorem~\ref{Thm:triangular_limit}.
\begin{enumerate}
\item[(1)] The fulfillment of this condition follows from the fact that $\inf_{t \in [0,1]}\big(\hat{H}_n(t) - Y_n(t)\big) \ge 0$ for all $n \in \NN$.
\item[(2)] Since $f_0$ is linear on $[0,1]$, $\hat{H}_n^{(2)}$ is convex on $[0,1]$ for every $n \in \NN$. The pointwise limit of convex functions is still convex, so $H^{(2)}_0$ is convex on $[1/m,1-1/m]$. The condition is then satisfied by letting $m \rightarrow \infty$.
\item[(3)] This condition always holds in view of our construction of $\hat{H}_n$ in (\ref{Eq:Def_H_n}).
\item[(4)] We consider two cases: (a) if $H^{(2)}_0(1-1/m) \rightarrow \infty$ as $m \rightarrow \infty$, then the conditions are satisfied in view of Corollary~\ref{Cor:triangular_limit}; (b) otherwise, it must be the case that $H^{(2)}_0(1^-)$ is bounded from above. Note that $H^{(2)}_0(1^-)$ is also bounded from below a.s., which can be proved by using Corollary~\ref{Cor:LS_triangular_uniform_rate_derivative} and the fact that $\hat{H}_n^{(2)}$ is convex. Denote by $\tau_{n_j}$ the knot of $\hat{f}_{n_j}$ closest to 1. Then by Lemma~2.2 of \citet{GJW2001b}, $\hat{H}_{n_j}(\tau_{n_j}) = Y_{n_j}(\tau_{n_j})$ and $\hat{H}_{n_j}'(\tau_{n_j}) = X_{n_j}(\tau_{n_j})$.  Since $t=1$ is a knot of $f_0$, consistency of $\hat{f}_{n_j}$ allows us to see that $\lim_{j \rightarrow \infty} \tau_{n_j} = 1$. Because $H_0^{(2)}(1^-)$ is finite and both $Y(t)$ and $X(t)$ are sample continuous processes, taking $\tau_{n_j} \rightarrow 1^-$ yields $H_0(1) = Y(1)$ and $H_0'(1) = X(1)$. Note that this argument remains valid even if $\tau_{n_j} > 1$, because in this scenario, $\hat{H}_{n_j}^{(2)}$ is linear and bounded on $[2-\tau_{n_j},\tau_{n_j}]$. 
\item[(5)] It follows from $\int_{1/m}^{1-1/m} \big(\hat{H}_{n_j}(t) - Y_{n_j}(t)\big)\, d\hat{H}_{n_j}^{(3)}(t) = 0$ that 
\[
\int_{1/m}^{1-1/m} \big(H_0(t) - Y(t)\big) \ dH_0^{(3)}(t) = 0.
\]
Since this holds for any $m$, one necessarily has that $\int_0^1 \big(H_0(t) - Y(t)\big) \ dH_0^{(3)}(t) = 0$.
\end{enumerate}

Consequently, in view of Theorem~\ref{Thm:triangular_limit}, the limit $Z_0$ is the same for any subsequences of $Z_n$ in $E_m$. Fix any $m > 1/\delta$. It follows that the full sequence $\{Z_n\}_n$ converges weakly in $E_m$ and has the limit $(H, H', H^{(2)}, H^{(3)}, Y, X)^T$. This, together with the fact that $H^{(3)}$ is continuous at any fixed $x_0 \in (0,1)$ with probability one (which can be proved using Conditions~(1) and (5) of $H$), yields (\ref{Eq:LS_pointwise_converg}).
\hfill $\Box$
\end{prooftitle}
\medskip

It can be inferred from Corollary~\ref{Cor:triangular_limit} and Theorem~\ref{Thm:LS_triangular_asymptotic} that both $\hat{f}_n(0)$ and $\hat{f}_n(1)$ do not converge to the truth at a $n^{-1/2}$ rate. In fact, \citet{Balabdaoui2007} proved that $\hat{f}_n(0)$ is an inconsistent estimator of $f_0(0)$. Nevertheless, the following proposition shows that $\hat{f}_n(0)$ is at most $O_p(1)$. For the case of the maximum likelihood estimator of a $k$ monotone density, we refer readers to \citet{GaoWellner2009} for a similar result.
\begin{prop}[Behavior at zero]
\label{Prop:LS_zero}
$\hat{f}_n(0) = O_p(1)$.
\end{prop}

\begin{prooftitle}{of Proposition~\ref{Prop:LS_zero}}
Let $\tau_n$ be the left-most point in $\mathcal{S}(\hat{f}_n)$. Since $\hat{f}_n(0)$ is finite for every $n$, the linearity of $\hat{f}_n$ on $[0,\tau_n]$ means that 
\[
\hat{f}_n(t) \ge (1-t/\tau_n)\hat{f}_n(0)
\]
for every $t \ge 0$. By Corollary~2.1 of \citet{GJW2001b},
\[
\mathbb{F}_n(\tau_n) = \int_0^{\tau_n} \hat{f}_n(t) dt \ge \int_0^{\tau_n} (1-t/\tau_n)\hat{f}_n(0) dt = \tau_n \hat{f}_n(0)/2.
\]
It follows from Theorem~9.1.2 of \citet{ShorackWellner1986} that
\[
\hat{f}_n(0) \le 2 \mathbb{F}_n(\tau_n)/\tau_n \le 2\sup_{t>0}\frac{\mathbb{F}_n(t)}{F_0(t)}\frac{F_0(t)}{t} \le O_p(1).
\]
\hfill $\Box$
\end{prooftitle}
\medskip

\subsubsection{Testing against a general convex density}
\label{Sec:test}
The asymptotic results established in Section~\ref{Sec:densityasym} can be used to test the linearity of a density function against a general convex alternative. To illustrate the main idea, in this section, we focus on the problem of testing $\mathcal{H}_0: f_0 (t) = 2(1-t) \mathbf{1}_{[0,1]}(t)$ against 
\[
\mathcal{H}_1: f_0 \in \mathcal{K} \; \mbox{ and } \; f_0 (t) \neq 2(1-t) \mathbf{1}_{[0,1]}(t) \mbox{ for some } t \in (0,\infty).
\]
The test we propose is free of tuning parameters. Since the triangular distribution is frequently used in practice, and is closely related to the uniform distribution (e.g. the minimum of two independent $U[0,1]$ is triangular), our test could be a valuable addition to practitioners' toolkit. For other tests on the linearity in the regression setting, we point readers to \citet{Meyer2003} and \citet{MeyerSen2013}.

Our test is based on the statistic 
\[
T_n := \sqrt{n}  \sup_{t \in [0,\infty)} \big\{2(1-t)\mathbf{1}_{[0,1]}(t) - \hat{f}_n(t) \big\}.
\]
The behavior of $T_n$ under $ \mathcal{H}_0$ is established below. 
\begin{prop}[Behavior of $T_n$ under $\mathcal{H}_0$]
\label{Prop:test_H_0}
Write $T:= -\inf_{t \in [0,1]} H^{(2)}(t)$, where $H$ is the invelope process defined in Theorem~\ref{Thm:triangular_limit}. Then $T_n \stackrel{d}{\rightarrow} T$. Moreover, $P(T \ge 0) = 1$.
\end{prop}
\begin{prooftitle}{of Proposition~\ref{Prop:test_H_0}}
The first part follows directly from Theorem~\ref{Thm:LS_triangular_asymptotic} and Corollary~\ref{Cor:triangular_limit}. For the second part, note that 
if $T<0$, then $\inf_{t \in [0,1]} H^{(2)}(t) > 0$, which would imply $H'(1) - H'(0) > 0$. But this would violate the characterization of $H$ (i.e. $H'(1) - H'(0) = X(1) - X(0) = 0$, since $X$ is a Brownian bridge). Consequently, $P(T \ge 0) = 1$.\hfill $\Box$
\end{prooftitle}
\medskip

\begin{figure}[ht!]
  \centering
  \includegraphics[scale=0.4,angle=90, trim = 180 0 150 0]{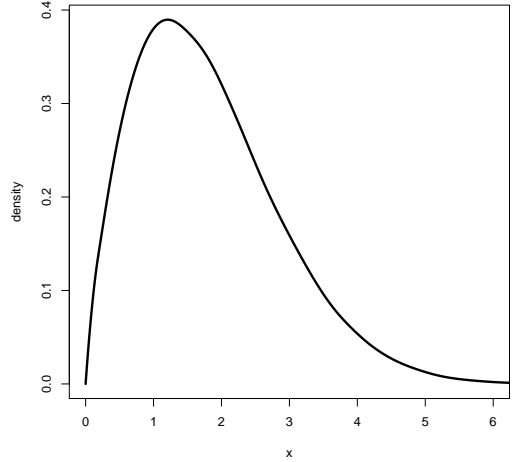} 
  \caption{The estimated density function of $T$ based on $10^5$ simulations. Here we used the kernel density estimator with the boundary correction at zero.}
\label{Fig:Test}
\end{figure}

The quantiles and the density function of $T$ can be approximated numerically using Monte Carlo methods. They are given in Table~\ref{Tab:Test_quantile} and Figure~\ref{Fig:Test}. Here we denote the upper $\alpha$ quantile of $T$ by $t_\alpha$. Our numerical evidence also suggests that the density function of $T$ appears to be log-concave. For related work on log-concavity of Chernoff's distribution (i.e. a different random variable also associated with a Brownian motion), see \citet{BalabdaouiWellner2014}.

\begin{table}[ht!]
  \centering
  \caption{Estimated upper-quantiles of $T$ based on $10^5$ simulations.}
  \label{Tab:Test_quantile}
  \begin{tabular}{| c | c  c  c  c  c|}\hline
      $\alpha$ & 1\% & 2.5\% & 5\% & 10\% & 20\% \\\hline
      $t_\alpha$ & 4.76 & 4.21 & 3.75 & 3.23 & 2.63\\\hline
  \end{tabular}
\end{table}

For a test of size $\alpha \in (0,1)$, we propose to reject $\mathcal{H}_0$ if $T_n > t_\alpha$. According to Proposition~\ref{Prop:test_H_0}, this test is asymptotically of size $\alpha$. In the following, we prove that our test is also consistent, i.e., its power goes to one as $n \rightarrow \infty$.

\begin{thm}[Consistency of the test]
\label{Thm:test_H_1}
Under $\mathcal{H}_1$, for any $\alpha > 0$, $P(T_n > t_\alpha) \rightarrow 1$ as $n \rightarrow \infty$.
\end{thm}

\begin{prooftitle}{of Theorem~\ref{Thm:test_H_1}}
Suppose that $f_0 (t) \neq 2(1-t) \mathbf{1}_{[0,1]}(t)$ for some $t \in (0,\infty)$. First, we show that there exists some $t^* \in (0,1)$ such that $f_0(t^*) < 2(1-t^*)$. Suppose the conclusion fails, then $f_0(t) \ge 2(1-t) \mathbf{1}_{[0,1]}(t)$ for all $t > 0$. Hence $\int_0^\infty f_0(t) dt \ge \int_0^\infty 2(1-t)\mathbf{1}_{[0,1]}(t) dt = 1$ with strict inequality if $f_0(t) > 2(1-t) \mathbf{1}_{[0,1]}(t)$ for some $t > 0$. But then $f_0$ is not a density and we conclude that $f_0(t) = 2(1-t) \mathbf{1}_{[0,1]}(t)$ for every $t \in (0,\infty)$. This contradiction yields the conclusion.

Next, it follows from Theorem~3.1 of \citet{GJW2001b} that the LSE $\hat{f}_n$ is consistent at $t^*$, i.e., $|\hat{f}_n(t^*) - f_0(t^*)|\stackrel{a.s.}{\rightarrow} 0$ . Therefore, almost surely, 
\[
T_n/\sqrt{n} >  2(1-t^*) - \hat{f}_n(t^*) =  2(1-t^*) - f_0(t^*) + f_0(t^*) - \hat{f}_n(t^*) \rightarrow 2(1-t^*) - f_0(t^*) > 0.
\]
It follows that $P(T_n > t_\alpha) \rightarrow 1$.
\hfill $\Box$
\end{prooftitle}
\medskip

\subsection{More general settings}
\label{Sec:densityext}
The aim of this subsection is to extend the conclusions presented in Section~\ref{Sec:densityspecial} to more general convex densities. 
We assume that $f_0$ is positive and linear on $(a,b)$ for some $0 \le a < b$, where the open interval $(a,b)$ is picked as the ``largest'' interval on which $f_0$ remains linear. More precisely, it means that there does not exist a bigger open interval $(a',b')$ (i.e. $(a,b) \subset (a',b') \subseteq (0,\infty)$) on which $f_0$ is linear.    

For the sake of notational convenience, we suppress the dependence of $H_*$ on $a$, $b$ and $F_0$ in the following two theorems. Their proofs are similar to those given in Section~\ref{Sec:densityspecial}, so are omitted for brevity. 
\begin{thm}[Characterization of the limit process]
\label{Thm:triangular_limit_extension}
Let $X(t) = \UU(F_0(t))$ and $Y(t) = \int_0^t X(s) ds$ for any $t > 0$. Then a.s., 
there exists a uniquely defined random continuously differentiable function $H_*$ on $[a,b]$ satisfying the following conditions:
\begin{enumerate}[(1)]
\item $H_*(t) \ge Y(t)$ for every $t \in [a,b]$;
\item $H_*$ has convex second derivative on $(a,b)$;
\item $H_*(a) = Y(a)$ and $H_*'(a) = X(a)$; 
\item $H_*(b) = Y(b)$ and $H_*'(b) = X(b)$; 
\item $\int_a^b \big(H_*(t)-Y(t)\big)\, dH_{*}^{(3)}(t) = 0$.
\end{enumerate}
\end{thm}

\begin{thm}[Rate and asymptotic distribution]
\label{Thm:triangular_dist_extension}
For any $0 < \delta \le (b-a)/2$, 
\[
	\sup_{x \in [a+\delta, b-\delta]} \Big(|\hat{f}_n(x) - f_0(x)|, \; |\hat{f}_n'(x) - f_0'(x)| \Big) = O_p(n^{-1/2}).
\]
Moreover,
\begin{align*}
\sqrt{n}
\left(\begin{array}{c}
\hat{f}_n (x) - f_0 (x)\\
\hat{f}_n' (x) - f_0' (x) 
\end{array} \right)
\Rightarrow
\left(\begin{array}{c}
H^{(2)}_*(x) \\
H^{(3)}_* (x)
\end{array} \right)  \mbox{ in }\, \mathcal{C}[a+\delta, b-\delta] \times \mathcal{D}[a+\delta, b-\delta],
\end{align*}
where $H_*$ is the invelope process defined in Theorem~\ref{Thm:triangular_limit_extension}.
\end{thm}

\subsection{Adaptation at the boundary points}
\label{Sec:densityadapt}
In this subsection, we study the pointwise convergence rate of the convex LSE $\hat{f}_n$ at the boundary points of the region where $f_0$ is linear. Examples of such points include $a$ and $b$ given in Section~\ref{Sec:densityext}. 

\begin{figure}[ht!]
  \centering
  \includegraphics[scale=0.35,angle=270]{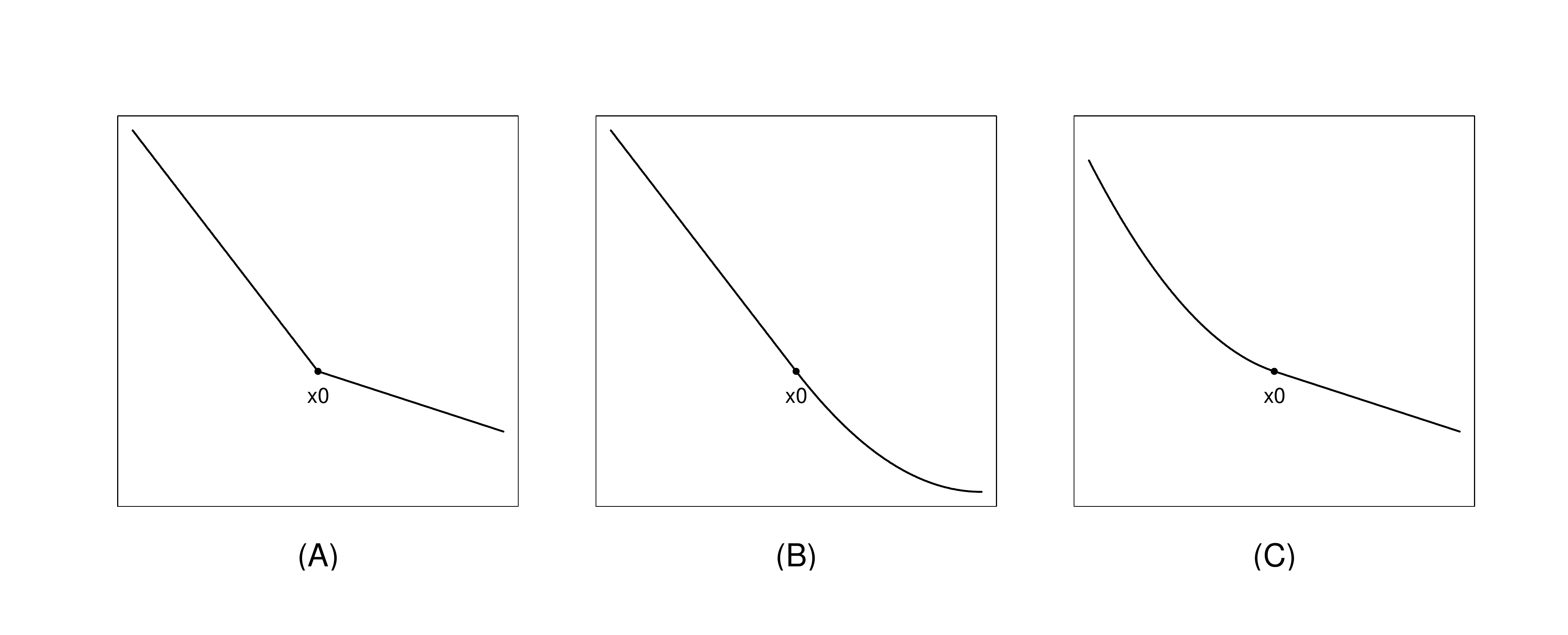} 
\caption{Typical $f_0$ in different cases are illustrated, where we set $\alpha = 2$ in (B) and (C).}
\label{fig:23ABC}
\end{figure}

To begin our discussion, we assume that $x_0 \in (0,\infty)$ is such a boundary point in the interior of the support (i.e. $f_0(x_0) > 0$). Here again $f_0$ is a convex (and decreasing) density function on $[0,\infty)$. Three cases are under investigation as below. These cases are illustrated in Figure~\ref{fig:23ABC}.

\begin{enumerate}[(A)]
\item $f_0(t) = f_0(x_0) + K_1 (t-x_0) + K_2 (t-x_0) \mathbf{1}_{[x_0,\infty)}(t)$ for every $t$ in a fixed (small) neighborhood of $x_0$, with $K_1 + K_2 < 0$ and $K_2 > 0$;
\item $f_0(t) = f_0(x_0) + K_1 (t-x_0) + K_2(t-x_0)^\alpha \mathbf{1}_{[x_0,\infty)}(t)$ for every $t$ in a fixed (small) neighborhood of $x_0$, with $K_1 < 0$, $K_2 > 0$ and $\alpha > 1$;
\item $f_0(t) = f_0(x_0) + K_1 (t-x_0) + K_2(x_0-t)^\alpha \mathbf{1}_{[0,x_0)}(t)$ for every $t$ in a fixed (small) neighborhood of $x_0$, with $K_1 < 0$, $K_2 > 0$ and $\alpha > 1$.
\end{enumerate}

Note that in all cases above, only the behavior of $f_0$ in a small neighborhood of $x_0$ is relevant. As pointed out in Example~2 of \citet{CaiLow2015}, the minimax optimal convergence rate at $x_0$ is $n^{-1/3}$ in (A). Furthermore, in (B) and (C), Example~4 of \citet{CaiLow2015} suggests that the optimal rate at $x_0$ is $n^{-\alpha/(2\alpha+1)}$. In the following, we prove that the convex LSE automatically adapts to optimal rates, up to a factor of $\sqrt{\log\log n}$. Though almost negligible, the factor of $\sqrt{\log\log n}$ here indicates that there might be room for improvement. These results should be viewed as a first step for the investigation of the adaptation of the convex LSE. One could also compare our results with \citet{Cator2011}, where the adaptation of the LSE in the context of decreasing density function estimation was tackled. 
\begin{thm}[Adaptation at the boundary points: I]
\label{Thm:LS_adapt_1}
In the case of (A), 
\[
|\hat{f}_n(x_0) - f_0(x_0)| = O_p\big(n^{-1/3}\sqrt{\log\log n}\big).
\]
Moreover,
\[
\min\big(\hat{f}_n(x_0) - f_0(x_0),0\big) = O_p(n^{-1/2}).
\]
\end{thm}

\begin{prooftitle}{of Theorem~\ref{Thm:LS_adapt_1}}
Suppose that for some fixed $\delta > 0$, (A) holds for every $t \in [x_0 - 2\delta, x_0 + 2\delta]$. Our first aim is to show that
\begin{align}
\label{Eq:LS_adapt_1_proof_1}
\inf_{t \in [x_0 - \delta, x_0 + \delta]} \min\big(\hat{f}_n(t) - f_0(t),0\big) = O_p(n^{-1/2}).
\end{align}

Suppose that $c:= \hat{f}_n(t) - f_0(t) < 0$. If we can find $x_0 -\delta \le t < x_0$ such that $\hat{f}_n$ and $f_0$ intersect at the point $(t,f_0(t))$, then 
\begin{align*}
	\big(\hat{F}_n(t) - \hat{F}_n(x_0 - \delta)\big)  -\big( F_0(t) - F_0(x_0 - \delta)\big) &\ge \frac{1}{2} \frac{|c|(t - x_0 + \delta)^2}{x_0-t} \ge 0; \\
	\big(\hat{F}_n(x_0) - \hat{F}_n(t)\big) - \big(F_0(x_0) - F_0(t)\big) &\le -\frac{1}{2}|c|(x_0-t) < 0.
\end{align*}
Otherwise, by taking $t=x_0 -\delta$, we can still verify the above two inequalities. Therefore, 
\begin{align*}
	O_p(n^{-1/2}) &= 2\lVert \hat{F}_n - F_0 \rVert_\infty \\
&\ge \max\Big\{\big|\big(\hat{F}_n(t) - F_0(t)\big) - \big(\hat{F}_n(x_0- \delta) - F_0(x_0- \delta)\big)\big|,\big|\big(\hat{F}_n(x_0) - F_0(x_0)\big) - \big(\hat{F}_n(t) - F_0(t)\big)\big|\Big\} \\
	&\ge \frac{|c|}{2}\max\bigg\{\frac{(t - x_0 + \delta)^2}{x_0-t},x_0-t\bigg\} \ge \frac{|c|\delta}{4},
\end{align*}
where the first equation follows from Marshall's lemma (see (\ref{Eq:marshall})). Consequently, (\ref{Eq:LS_adapt_1_proof_1}) holds true. As a remark, the above arguments are essentially the same as those illustrated in the proof of Theorem~\ref{Thm:LS_triangular_rate}, where one only relies on the fact that $f_0$ is linear on $[x_0 - \delta, x_0]$ (i.e. ``one-sided linearity'').

In the rest of the proof, it suffices to only consider the situation of $\hat{f}_n(x_0) - f_0(x_0) > 0$. Recall that to handle this scenario, the proof of Theorem~\ref{Thm:LS_triangular_rate} makes use of the fact that the triangular density is linear on the whole $[0,1]$ (i.e. ``two-sided linearity''). This is no longer true here since $f_0$ is not linear on $[x_0-\delta,x_0+\delta]$. In the following, we deploy a different strategy to establish the rate. 

Let $\tau_n^-  = \max\{t \in \mathcal{S}(\hat{f}_n), t < x_0\}$ and $\tau_n^+ = \min\{t \in \mathcal{S}(\hat{f}_n), t \ge x_0\}$, where $\mathcal{S}(\cdot)$ is defined in (\ref{Eq:defknots}). We consider three different cases separately.
\begin{enumerate}[(a)]
\item $\max\big(\tau_n^+ - x_0,x_0 - \tau_n^-\big) > \delta$. Then $f_0$ and $\hat{f}_n$ are linear on either $[x_0 - \delta, x_0]$ or $[x_0, x_0+\delta]$. It follows from the line of reasoning as in the proof of Theorem~\ref{Thm:LS_triangular_rate} that $\max\big(\hat{f}_n(x_0) - f_0(x_0),0\big) = O_p(n^{-1/2})$.

\item $\max\big(\tau_n^+ - x_0,x_0 - \tau_n^-\big) \le \delta$ and $\tau_n^+ - \tau_n^- > n^{-1/3}$. Note that being in the set $\mathcal{S}(\hat{f}_n)$ implies that $\mathbb{F}_n(\tau_n^-) = \hat{F}_n(\tau_n^-)$ and $\mathbb{F}_n(\tau_n^+) = \hat{F}_n(\tau_n^+)$ for $\tau_n^-$ and $\tau_n^+$. Since $\hat{f}_n$ is linear on $[\tau_n^-,\tau_n^+]$,
\begin{align}
\label{Eq:LS_adapt_1_proof_2} &\Big(\mathbb{F}_n(\tau_n^+) - F_0(\tau_n^+)\Big) - \Big(\mathbb{F}_n(\tau_n^-) - F_0(\tau_n^-)\Big) \\
\notag&= \int_{\tau_n^-}^{\tau_n^+} \big(\hat{f}_n(t) - f_0(t)\big)dt \\
\notag&\ge \big(\hat{f}_n(x_0) - f_0(x_0)\big)(\tau_n^+ - \tau_n^-)/2 + \min\Big(\hat{f}_n(\tau_n^-) - f_0(\tau_n^-),\hat{f}_n(\tau_n^+) - f_0(\tau_n^+),0\Big)(\tau_n^+ - \tau_n^-).
\end{align}
By the law of the iterated logarithm for local empirical processes (cf. Lemma 4.3.5 of \citet{CsorgoHorvath1993}), (\ref{Eq:LS_adapt_1_proof_2}) is at most $\sqrt{\tau_n^+ - \tau_n^-} O_p\big(n^{-1/2} \sqrt{\log\log n}\big)$. In view of (\ref{Eq:LS_adapt_1_proof_1}), rearranging the terms in the above inequality yields
\[
\hat{f}_n(x_0) - f_0(x_0) \le (\tau_n^+ - \tau_n^-)^{-1/2} O_p\big(n^{-1/2}\sqrt{\log\log n}\big) + O_p(n^{-1/2}) \le O_p\big(n^{-1/3}\sqrt{\log\log n}\big).
\]

\item $\tau_n^+ - \tau_n^- \le n^{-1/3}$. Now we define $\tau_n^{++} = \max\{t \in \mathcal{S}(\hat{f}_n), t <  x_0 + 2 n^{-1/3}\}$ and $\tau_n^{+++} = \min\{t \in \mathcal{S}(\hat{f}_n), t \ge x_0 + 2 n^{-1/3}\}$. Here the existence of $\tau_n^{++}$ is guaranteed by the condition that $\tau_n^+ - \tau_n^- \le n^{-1/3}$. Furthermore, we note that in our definitions $\tau_n^{+}$ and $\tau_n^{++}$ might not be distinct. Within this setting, three further scenarios are to be dealt with.
\begin{enumerate}[(c1)]
\item $\tau_n^{+++} - x_0 > \delta$. Since both $f_0$ and $\hat{f}_n$ are linear on $[\tau_n^{++},\tau_n^{+++}]$, one can apply a strategy similar to that used in the proof of Theorem~\ref{Thm:LS_triangular_rate} to show that $|\hat{f}_n(\tau_n^{++}) - f_0(\tau_n^{++})| = O_p(n^{-1/2})$. Consequently,
\[
\hat{f}_n(x_0) - f_0(x_0) = \hat{f}_n(\tau_n^{++}) - f_0(\tau_n^{++}) - \int_{x_0}^{\tau_n^{++}} \big(\hat{f}_n'(t) - f_0'(t)\big)dt = O_p(n^{-1/3}),
\]
where we have used the facts that $|\tau_n^{++} - x_0| < 2n^{-1/3}$ and $\sup_{t \in [x_0,x_0+\delta]}|\hat{f}_n'(t) - f_0'(t)| = O_p(1)$. Here the second fact can be derived by invoking
\[
K_1 + o_p(1) \le \inf_{t \in [x_0 - \delta, x_0 + \delta]}\hat{f}_n'(t) \le \sup_{t \in [x_0 - \delta, x_0 + \delta]}\hat{f}_n'(t) \le K_1 +K_2 + o_p(1),
\]
which follows easily from consistency of $\hat{f}_n'$ in estimating $f_0'$ at the points $x_0 \pm \delta$.
\item $3n^{-1/3} < \tau_n^{+++} - x_0 \le \delta$. Then $|\tau_n^{+++} - \tau_n^{++}| > n^{-1/3}$. Using essentially the same argument as in (b), we see that
\[
\max\big(\hat{f}_n(\tau_n^{++}) - f_0(\tau_n^{++}), 0\big) =  O_p\big(n^{-1/3}\sqrt{\log\log n}\big),
\]
and hence,
\[
|\hat{f}_n(\tau_n^{++}) - f_0(\tau_n^{++})| = O_p\big(n^{-1/3}\sqrt{\log\log n}\big).
\]
We then apply the argument presented in (c1) to derive
\[
|\hat{f}_n(x_0) - f_0(x_0)| = O_p\big(n^{-1/3}\sqrt{\log\log n}\big).
\]
\item $\tau_n^{+++} - x_0 \le 3n^{-1/3}$. Then $n^{-1/3} \le \tau_n^{+++} - \tau_n^{+} \le 2n^{-1/3}$. By proceeding as in the proof of Lemma~4.3 of \citet{GJW2001b}, we are able to verify that
\[
\inf_{t \in [\tau_n^{+}, \tau_n^{+++}]} |\hat{f}_n(t) - f_0(t)| = O_p(n^{-1/3}).
\]
Here we invoked Lemma~A.1 of \citet{BalabdaouiWellner2007} with $k = 1$ and $d = 2$ to verify the above claim. See also \citet{KimPollard1990}. Finally, we can argue as in (c1) to show that $\hat{f}_n(x_0) - f_0(x_0) = O_p(n^{-1/3})$.
\end{enumerate}
\end{enumerate}
The proof is complete by taking into account all the above cases. We remark that in our proof the worst case scenarios that drive the convergence rate are (b) and (c2).

\hfill $\Box$
\end{prooftitle}
\medskip

\begin{thm}[Adaptation at the boundary points: II]
\label{Thm:LS_adapt_2}
In the case of (B) or (C), 
\[
|\hat{f}_n(x_0) - f_0(x_0)| = O_p\big(n^{-\alpha/(2\alpha+1)}\sqrt{\log\log n}\big).
\]
Moreover,
\[
\min\big(\hat{f}_n(x_0) - f_0(x_0),0\big) = O_p(n^{-1/2}).
\]
\end{thm}
The proof of Theorem~\ref{Thm:LS_adapt_2} is of more technical nature and is deferred to Appendix II.

\section{Estimation of a regression function}
\label{Sec:regress}
Changing notation slightly from the previous section, we now assume that we are 
given pairs $\{(X_{n,i}, Y_{n,i}): i = 1,\ldots,n\}$ for $n = 1,2,\ldots$ with
\[
Y_{n,i} = r_0(X_{n,i}) + \epsilon_{n,i}, \quad \mbox{ for } i = 1,\ldots,n,
\]
where $r_0: [0,1] \rightarrow \RR$ is a convex function, and where $\{\epsilon_{n,i}:i=1,\ldots,n\}$ is a triangular array of IID random variables satisfying 
\begin{enumerate}[(i)]
\item $\mathbb{E} \epsilon_{1,1} = 0$ and $\sigma_0^2 = \mathbb{E} \epsilon_{1,1}^2 < \infty$.
\end{enumerate}
To simplify our analysis, the following fixed design is considered:
\begin{enumerate}[(i)]
\setcounter{enumi}{1}

\item $X_{n,i} = i/(n+1)$ for $i = 1,\ldots,n$.
\end{enumerate}

The LSE of $r_0$ proposed by \citet{BalabdaouiRufibach2008} is
\begin{align}
\label{Eq:RE_estimator}
\hat{r}_n = \argmin_{g \in \mathcal{K}} \left( \frac{1}{2}\int_0^1 g(t)^2 dt - \frac{1}{n}\sum_{i=1}^n g(X_{n,i}) Y_{n,i} \right),
\end{align}
where, in this section, $\mathcal{K}$ denotes the set of all continuous convex functions on $[0,1]$.

We note that the above estimator is slightly different from the more ``classical'' LSE that
minimizes $\frac{1}{2n} \sum_{i=1}^n g(X_{n,i})^2 - \frac{1}{n}\sum_{i=1}^n g(X_{n,i}) Y_{n,i}$ over $\mathcal{K}$.
Since its criterion function has a completely discrete nature, different techniques are needed to prove analogous results.
We will not pursue this direction in the manuscript. 

\subsection{Basic properties}
In this subsection, we list some basic properties of $\hat{r}_n$ given as (\ref{Eq:RE_estimator}).
For any $t \in [0,1]$, define
\begin{align*}
R_0(t) &= \int_0^t r_0(s) ds; \\
\hat{R}_n(t) &= \int_0^t \hat{r}_n(s) ds; \\
\mathbb{R}_n(t) &= \frac{1}{n}\sum_{i=1}^n Y_{n,i} \mathbf{1}_{[X_{n,i},\, \infty)}(t). \\
\end{align*}

\begin{prop}
\label{Prop:RE_exist}
The LSE $\hat{r}_n$ exists and is unique.
\end{prop}

\begin{prop}
\label{Prop:RE_prop}
Let $\mathcal{S}(\hat{r}_n)$ denote the set of knots of $\hat{r}_n$. The following properties hold:
\begin{enumerate}[(a)]
\item
$\int_0^t \big(\hat{R}_n(s) -\mathbb{R}_n(s)\big)\, ds \left \{ \begin{array}{l}{= 0, \; \mbox{ for all } t \in \mathcal{S}(\hat{r}_n)\cup \{0, 1\},}\\{\ge 0, \;  \mbox{ for all } t \in [0,1];}\end{array}\right.$ 
\item $\hat{R}_n(t) = \mathbb{R}_n(t)$ for any $t \in \mathcal{S}(\hat{r}_n)\cup \{0,1\}$;
\item $\int_0^1 \big\{\int_0^t \big(\hat{R}_n(s) -\mathbb{R}_n(s)\big) ds \big\} \, d\hat{r}_n'(t) = 0$.
\end{enumerate}
\end{prop}

\begin{prop}[Marshall's lemma]
\label{Prop:RE_marshall}
\begin{align}
\label{Eq:RE_marshall}
	\lVert \hat{R}_n - R_0 \rVert_{\infty} \le 2 \lVert \mathbb{R}_n - R_0 \rVert_{\infty}.
\end{align}
\end{prop}

\begin{prooftitle}{of Propositions~\ref{Prop:RE_exist}, \ref{Prop:RE_prop} and \ref{Prop:RE_marshall}}
Note that we can rewrite $\hat{r}_n$ as 
\[
\argmin_{g \in \mathcal{K}} \Big( \frac{1}{2}\int_0^1 g(t)^2 dt - \int_0^1 g(t) d \RR_n(t) \Big).
\]
From this perspective, it is easy to check that Proposition~\ref{Prop:RE_exist}, Proposition~\ref{Prop:RE_prop} and Proposition~\ref{Prop:RE_marshall} follow from, respectively, slight modifications of Lemma~2.1, Lemma~2.2 of \citet{GJW2001b} and Theorem~1 of \citet{DRW2007}.
\hfill $\Box$
\end{prooftitle}

We remark that this version of Marshall's lemma in the regression setting serves as an important tool to establish consistency and the rate.
In particular, in the following, we show that (\ref{Eq:RE_marshall}) easily yields consistency of $\hat{R}_n$, and consistency of $\hat{r}_n$ follow from this together with convexity of $\hat{r}_n$.

\begin{prop}[Consistency]
\label{Prop:RE_consistency}
Under Assumptions (i) -- (ii), for any $0 < \delta \le 1/2$, 
\[
	\sup_{x \in [\delta, 1-\delta]} |\hat{r}_n(x) - r_0(x) |  \stackrel{\mbox{a.s.}}{\rightarrow} 0, \quad\mbox{ as } n \rightarrow \infty.
\]
\end{prop}

\begin{prooftitle}{of Proposition~\ref{Prop:RE_consistency}}
The proof replies on the following intermediate result: suppose there are differentiable functions $G_0, G_1,\ldots$, all are $[0,1]\rightarrow \RR$ and have convex first derivative, then for any fixed $0 < \delta \le 1/2$, $\lVert G_k - G_0 \rVert_{\infty} \rightarrow 0$ implies $\sup_{t \in [\delta, 1-\delta]} |G_k'(t) - G_0'(t)| \rightarrow 0$.

To see this, we first show that 
\begin{align}
\label{Eq:RE_marshall1}
\limsup_{k \rightarrow \infty} \sup_{t \in [\delta, 1-\delta]} \{ G_k'(t) - G_0'(t) \} \le 0. 
\end{align}
Note that for any $t \in [\delta,1-\delta]$ and any $0 < \epsilon \le \delta$,
\begin{align*}
G_k'(t) - G_0'(t) \le \max \bigg( \frac{G_k(t+\epsilon) - G_k(t)}{\epsilon} - G_0'(t), \frac{G_k(t) - G_k(t-\epsilon)}{\epsilon} - G_0'(t) \bigg).
\end{align*}
It then follows that
\[
\limsup_{k \rightarrow \infty} \sup_{t \in [\delta, 1-\delta]} \{ G_k'(t) - G_0'(t) \} \le \sup_{t \in [\delta, 1-\delta]} \bigg( \frac{G_0(t+\epsilon) - G_0(t)}{\epsilon} - G_0'(t), \frac{G_0(t) - G_0(t-\epsilon)}{\epsilon} - G_0'(t) \bigg).
\]
Our claim (\ref{Eq:RE_marshall1}) can be verified by letting $\epsilon \rightarrow 0$.

Secondly, we show that 
\begin{align*}
\liminf_{k \rightarrow \infty} \inf_{t \in [0, 1]} \{ G_k'(t) - G_0'(t) \} \ge 0. 
\end{align*}
We prove this by contradiction. Suppose that $\liminf_{k \rightarrow \infty} \inf_{t \in [0, 1]} \{ G_k'(t) - G_0'(t) \} = -M$ for some $M > 0$. By extracting subsequences if necessary, we can assume that $\inf_{t \in [0, 1]} \{ G_k'(t) - G_0'(t) \} \rightarrow -M$ as $k \rightarrow \infty$. In view of (\ref{Eq:RE_marshall1}), it follows from the convexity of $G_k'$ and $G_0'$ that one can find an interval $\mathcal{I}_k$ of positive length $\delta$ (which can depend on $M$) such that $\inf_{t \in \mathcal{I}_k} \{ G_k'(t) - G_0'(t) \} \le -M/2$ for every $k > K$, where $K$ is a sufficiently large integer. This implies that
\[
\lVert G_k - G_0 \rVert_{\infty} \ge \sup_{t \in \mathcal{I}_k} |G_k(t) - G_0(t)| \ge M\delta/4 > 0,
\]
for every $k > K$, which contradicts the fact that $\lVert G_k - G_0 \rVert_{\infty} \rightarrow 0$ as $k \rightarrow \infty$. 

Combining these two parts together completes the proof of the intermediate result.

Since $\lVert \mathbb{R}_n - R_0 \rVert_{\infty} \stackrel{a.s.}{\rightarrow} 0$ by empirical process theory, Proposition~\ref{Prop:RE_marshall} entails that $\lVert \hat{R}_n - R_0 \rVert_{\infty} \stackrel{a.s.}{\rightarrow} 0$. Consistency of $\hat{r}_n$ then follows straightforwardly from the above intermediate result.
\hfill $\Box$
\end{prooftitle}

\subsection{Rate of convergence and asymptotic distribution}
In the following, we assume that $r_0$  that is linear on $(a,b) \subseteq (0,1)$. 
Moreover, $(a,b)$ is ``largest'' in the sense that one can not find a bigger open interval $(a',b')$ on which $r_0$ remains linear. 
\begin{thm}[Rate and asymptotic distribution]
\label{Thm:RE_rate_dist}
Under Assumptions (i) -- (ii), for any $0 < \delta \le (b-a)/2$, 
\[
	\sup_{x \in [a+\delta, b-\delta]} \Big(|\hat{r}_n(x) - r_0(x)|, \; |\hat{r}_n'(x) - r_0'(x)| \Big) = O_p(n^{-1/2}).
\]
Moreover,
\begin{align*}
\sqrt{n}
\left(\begin{array}{c}
\hat{r}_n (x) - r_0 (x)\\
\hat{r}_n' (x) - r_0' (x) 
\end{array} \right)
\Rightarrow
\sigma_0
\left(\begin{array}{c}
(b-a)^{1/2} \, \tilde{H}^{(2)}\big(\frac{x-a}{b-a}\big) \\
(b-a)^{-1/2} \, \tilde{H}^{(3)}\big(\frac{x-a}{b-a}\big) 
\end{array} \right)  \mbox{ in }\, \mathcal{C}[a+\delta, b-\delta] \times \mathcal{D}[a+\delta, b-\delta],
\end{align*}
where $\tilde{H}$ is the invelope process defined in the second part of Theorem~\ref{Thm:triangular_limit}. 
\end{thm}
The proof of Theorem~\ref{Thm:RE_rate_dist} is very similar to what has already been shown in Section~\ref{Sec:densityspecial}, so is omitted for the sake of brevity.

In presence of the linearity of $r_0$ on $(a,b)$, the limit distribution of the process $\sqrt{n}(\hat{r}_n - r_0)$ on $(a,b)$ does not depend on $r_0(a)$ or $r_0(b)$. 
In addition, the above theorem continues to hold if we weaken Assumption (ii) to:
\begin{enumerate}[(i')]
\setcounter{enumi}{1}
\item $\sup_{t \in [0,1]}\Big|\frac{1}{n}\sum_{i=1}^n \mathbf{1}_{[X_{n,i},\, \infty)}(t) - t \Big| = o\big(n^{-1/2}\big)$.
\end{enumerate}

Theoretical results in the random design are also possible, where for instance, we can assume that $\{X_{n,i}, i = 1,\ldots,n\}$ are IID uniform random variables on $[0,1]$. In this case, Theorem~\ref{Thm:RE_rate_dist} is still valid, while a process different from $\tilde{H}$ is required to characterize the limit distribution. This follows from the fact that in the random design $\sqrt{n}(\mathbb{R}_n - R_0)$ can converge to a Gaussian process that is not a Brownian motion.

\section{Appendices: proofs}
\subsection{Appendix I: existence of the limit process}
Recall that $f_0(t) = 2(1-t)\mathbf{1}_{[0,1]}(t)$, $X(t) = \mathbb{U}(F_0(t))$, and $Y(t) = \int_0^t X(s)ds$ for $t \ge 0$. In this section, we show the existence and uniqueness of the invelope process $H$. The case of $\tilde{H}$ can be handled using essentially the same arguments, so is omitted here.
 
Lemmas~\ref{Lem:limsolution} -- \ref{Lem:limequicontinuous} are needed to prove Theorem~\ref{Thm:triangular_limit}.
\begin{lem}
\label{Lem:limsolution}
Let the functional $\phi(g)$ be defined as
\begin{align*}
	\phi(g) = \frac{1}{2}\int_{0}^{1}g^2(t) \, dt - \int_{0}^{1}g(t) \, dX(t)
\end{align*}
for functions in the set $\mathcal{G}_k = \{ g:[0,1]\rightarrow \RR, g \mbox{ is convex}, g(0) = g(1) = k \}$. 
Then with probability one, the problem of minimizing $\phi(g)$ over $\mathcal{G}_k$ has a unique solution.
\end{lem}
\begin{prooftitle}{of Lemma~\ref{Lem:limsolution}}
We consider this optimization problem in the metric space $L^2$. First, we show that if it exists, the minimizer must be in the subset
\begin{align*}
	\mathcal{G}_{k,M} = \{ g:[0,1]\rightarrow \RR, g \mbox{ is convex}, g(0) = g(1) = k, \inf_{[0,1]} g(t) \ge - M \}
\end{align*}
for some $0 < M < \infty$. To verify this, we need the following result
\begin{align}
\label{Eq:GCfinite}
	\sup_{g \in \mathcal{G}_{1,1}}\left|\int_0^1 g(t) \, dX(t)\right| < \infty,\quad \mbox{a.s.} 
\end{align}
Let $W(t)$ be a standard Brownian motion. We note that $\int_0^1 g(t) \, dX(t)$ has the same distribution as
\[
\int_0^1 g\big(1-\sqrt{1-t}\big) dW(t) - W(1) \int_0^1 g(t) f_0(t) dt.
\]
Using the entropy bound of $\mathcal{G}_{1,1}$ in $L^2$ (Theorem 2.7.1 of \citet{GuntuboyinaSen2013}) and Dudley's theorem (cf. Theorem~2.6.1 of \citet{Dudley1999}), we can establish that 
\[
\mathcal{H} = \big\{h:[0,1]\rightarrow \RR \; | \; h(t) = g\big(1-\sqrt{1-t}\big), g \in \mathcal{G}_{1,1} \big\}
\]
is a GC-set. As $\int_0^1 g\big(1-\sqrt{1-t}\big) dW(t)$ is an isonormal Gaussian process indexed by $\mathcal{H}$, we have that a.s. $\sup_{g \in \mathcal{G}_{1,1}}\big|\int_0^1 g\big(1-\sqrt{1-t}\big) dW(t)\big| < \infty$. Furthermore, it is easy to check that $W(1) < \infty$ a.s. and $\sup_{g \in \mathcal{G}_{1,1}}\big| \int_0^1 g(t) f_0(t)dt\big| \le 2$. So our claim of (\ref{Eq:GCfinite}) holds. 

Now for sufficiently large $M$ (with $M > k$), 
\begin{align*}
	\sup_{g \in \mathcal{G}_{k,M}}\left|\int_{0}^{1}g(t) \, dX(t)\right| \le M \sup_{g \in \mathcal{G}_{1,1}}\left|\int_0^1 g(t) \, dX(t)\right|.
\end{align*}
Thus, for any $g \in \mathcal{G}_k$ with $\inf_{[0,1]} g = -M$, $\int_{0}^{1}g(t) dX(t)$ is at most $O(M)$. Furthermore, $\frac{1}{2}\int_{0}^{1}g^2(t) \, dt$ is of $O(M^2)$. Since $\phi$ at the minimizer could at most be as large as $\phi(0) = 0$, we conclude that it suffices to only consider functions in $\mathcal{G}_{k,M}$ for some sufficiently large $M$.

Note that the functional $\phi$ is continuous (cf. Dudley's theorem) and strictly convex. Moreover, for $g_1, g_2 \in \mathcal{G}_{k, M}$, if $\int_0^1 \big(g_1(t) - g_2(t)\big)^2 dt = 0$, then $g_1 = g_2$ on $[0,1]$. Since $\mathcal{G}_{k,M}$ is compact in $L^2$, the existence and uniqueness follow from a standard convex analysis argument in a Hilbert space. 
\hfill $\Box$
\end{prooftitle}
\medskip

As a remark, it can be seen from the proof of Lemma~\ref{Lem:limsolution} that for a given $\omega \in \Omega$ from the sample space (which determines the value of $X(t)$), if the function $\phi$ has a unique minimizer over $\mathcal{G}_1$ (which happens a.s.), it also admits a unique minimizer over $\mathcal{G}_k$ for any $k > 1$.

\begin{lem}
\label{Lem:parabolictangent}
Almost surely, $Y(t)$ does not have parabolic tangents at either $t = 0$ or $t = 1$.
\end{lem}
\begin{prooftitle}{of Lemma~\ref{Lem:parabolictangent}}
First, consider the case of $t=0$. Theorem~1 of \citet{Lachal1997} says that
\[
\limsup_{t \rightarrow 0^+} \frac{\big|\int_0^t W(s) ds\big|}{\sqrt{(2/3)t^3\log\log(1/t)}} = 1, \quad a.s.,
\]
where $W$ is a standard Brownian motion. From this, it follows that
\begin{align*}
\limsup_{t \rightarrow 0^+} \frac{|Y(t)|}{\sqrt{(2/3)t^3\log\log(1/t)}} &= \limsup_{t \rightarrow 0^+} \frac{\big|\int_0^t \{W(2s-s^2) - (2s-s^2)W(1)\}ds\big|}{\sqrt{(2/3)t^3\log\log(1/t)}} \\
&= \limsup_{t \rightarrow 0^+} \frac{\big|\int_0^t \sqrt{2}W(s)ds\big|}{\sqrt{(2/3)t^3\log\log(1/t)}} = \sqrt{2}, \quad a.s.
\end{align*}
thanks to the scaling properties of Brownian motion ($W(at) \stackrel{d}{=} \sqrt{a} W(t)$). This implies that $Y(t)$ does not have a parabolic tangent at $t=0$.

Second, consider the case of $t=1$. Note that  $\lim_{t\rightarrow 1^-} \frac{Y(1) - Y(t)}{1-t} = X(1) = 0$ and 
\begin{align*}
\limsup_{t\rightarrow 1^-} \bigg| \frac{Y(1) - Y(t)}{(1-t)^2}\bigg| & = \limsup_{t\rightarrow 1^-}  \bigg| \frac{\int_{t}^1 \big\{W(1)-W(F_0(s))\big\}ds - W(1)(1-t)^3/3 }{(1-t)^2}\bigg|\\
& = \lim_{t\rightarrow 1^-} \bigg|\frac{ \int_{t}^1 \big\{W(1)-W(F_0(t))\big\}ds }{(1-t)^2}\bigg| \stackrel{a.s.}{=} \lim_{t\rightarrow 0^+} \bigg|\frac{\int_{0}^t W(s^2)ds}{t^2}\bigg|,
\end{align*}
where $F_0(t) =  2t - t^2$. Therefore, to prove that $Y(t)$ does not have a parabolic tangent at $t=1$, it suffices to show that
\[
\limsup_{t\rightarrow 0^+}\bigg| \frac{\int_0^t W(s^2)ds}{t^2} \bigg| = \infty, \quad a.s.
\]
Denote by $Z(t) = \int_0^t W(s^2)ds$. For any $0 < t_1 < t_2 < 1$, we argue that the random variable
\[
Z(t_2)-Z(t_1) - (t_2-t_1)W(t_1^2)
\]
follows a distribution of $N\Big(0, t_2^4 / 6 - t_1^2 t_2^2 + 4t_1^3 t_2 / 3 - t_1^4 / 2\Big)$.
This is because
\begin{align*}
Z(t_2)-Z(t_1) - (t_2-t_1)W(t_1^2) &= \int_{t_1}^{t_2} \{W(s^2) - W(t_1^2)\} ds \\
& = \int_{t_1^2}^{t_2^2} \frac{W(s) - W(t_1^2)}{2\sqrt{s}} ds = \int_{t_1^2}^{t_2^2} \frac{1}{2\sqrt{s}} \int_{t_1^2}^s dW(u) ds \\
& = \int_{t_1^2}^{t_2^2} \int_u^{t_2^2} \frac{1}{2\sqrt{s}} ds\; dW(u) =  \int_{t_1^2}^{t_2^2} (t_2-\sqrt{u}) dW(u),
\end{align*}
where we invoked the stochastic Fubini's theorem in the last line, and thus,
\[
\mathrm{Var}\Big\{Z(t_2)-Z(t_1) - (t_2-t_1)W(t_1^2) \Big\} = \int_{t_1^2}^{t_2^2} (t_2-\sqrt{u})^2 du = t_2^4 / 6 - t_1^2 t_2^2 + 4t_1^3 t_2 / 3 - t_1^4 / 2. 
\]

Now setting $t_1 = 1/2$ and $t_{i+1} = t_i^2$ for every $i \in \mathbb{N}$. It is easy to check that the collection of random variables
\[
\Big\{Z(t_{i-1})-Z(t_i) - (t_{i-1}-t_i)W(t_i^2)\Big\}_{i=1}^\infty
\]
is mutually independent, so  
\begin{align}
\label{Eq:para_tangent}
\limsup_{i \rightarrow \infty} \frac{Z(t_{i-1})-Z(t_i) - (t_{i-1}-t_i)W(t_i^2)}{t_{i-1}^2} = \infty, \quad a.s.,
\end{align}
where we made use of the fact that 
\[
\lim_{i \rightarrow \infty} \frac{\sqrt{t_{i-1}^4/6 - t_{i}^2t_{i-1}^2 + 4t_i^3 t_{i-1}/3 - t_i^4/2}}{t_{i-1}^2} = 1/\sqrt{6}.
\]

Assume that there exists some $K > 0$ such that $|Z(t)| \le K t^2$ for all sufficiently small $t>0$. But it follows from (\ref{Eq:para_tangent}) that a.s. one can find a subsequence of $\mathbb{N}$ (denoted by $\{i_j\}_{j=1}^\infty$) satisfying
\[
\frac{Z(t_{i_j-1})-Z(t_{i_j}) - (t_{i_j-1}-t_{i_j})W(t_{i_j}^2)}{t_{i_j-1}^2}  \ge K+1, \quad \mbox{for } j = 1,2,\ldots
\]
Consequently,
\begin{align*}
\frac{Z(t_{i_j-1})}{t_{{i_j}-1}^2} &\ge K+1 + \frac{Z(t_{i_j})}{t_{{i_j}-1}^2} + \frac{W(t_{i_j}^2)(t_{{i_j}-1}-t_{i_j})}{t_{{i_j}-1}^2} \\
& \ge K+1 - K t_{i_j} - \frac{|W(t_{i_j}^2)|}{\sqrt{t_{i_j}}} \frac{\sqrt{t_{i_j}}(t_{{i_j}-1}-t_{i_j})}{t_{i_j}} \rightarrow K+1, 
\end{align*}
as $j \rightarrow \infty$. The last step is due to the facts of $t_{i_j} \rightarrow 0^+$ and $\limsup_{s\rightarrow 0^+}|W(s)|/s^{1/4} = 0$ a.s. (which is a direct application of the law of the iterated logarithm). The proof is completed by contradiction.
\hfill $\Box$
\end{prooftitle}
\medskip

Now denote by $f_k$ the unique function which minimizes $\phi(g)$ over $\mathcal{G}_k$. Let $H_k$ be the second order integral satisfying $H_k(0) = Y(0)$, $H_k(1) = Y(1)$ and $H_k^{(2)} = f_k$. 
\begin{lem}
\label{Lem:limprop}
Almost surely, for every $k \in \mathbb{N}$, $f_k$ and $H_k$ has the following properties:
\begin{enumerate}[(i)]
\item $H_k(t) \ge Y(t)$ for every $t \in [0,1]$;
\item $\int_0^1 \big(H_k(t)-Y(t)\big)\, df_k'(t) = 0$, where the derivative can be interpreted as either the left or the right derivative;
\item $H_k(t) = Y(t)$ and $H_k'(t) = X(t)$ for any $t \in \mathcal{S}(f_k)$, where $\mathcal{S}$ is the set of knots;
\item $\int_0^t f_k(s) ds \le X(t) - X(0)$ and $\int_t^1 f_k(s) ds \le X(1) - X(t)$ for any $t \in \mathcal{S}(f_k)$;
\item $f_k$ is a continuous function on $[0,1]$.
\end{enumerate}
\end{lem}
\begin{prooftitle}{of Lemma~\ref{Lem:limprop}}
To show (i), (ii), (iii) and (iv), one may refer to Lemma~2.2 and Corollary~2.1 of \citet{GJW2001a} and use a similar functional derivative argument. 

For (v), we note that since $f_k$ is convex, discontinuity can only happen at $t=0$ or $t=1$. In the following, we show that it is impossible at $t=0$. Suppose that $f_k$ is discontinuous at zero. Consider the class of functions $g_{\delta}(t) = \max(1-t/\delta,0)$. Then 
\[
\Big(f_k(t) + \epsilon g_{\delta}(t)\Big)\mathbf{1}_{(0,1)}(t) + k \mathbf{1}_{\{0,1\}}(t) \in \mathcal{G}_k
\] 
for every $\epsilon \in (0, k - \liminf_{t \rightarrow 0^+}f_k(t)]$. By considering the functional derivative of $\phi(g)$ and using integration by parts, we obtain that for any $\delta >0$,
\begin{align*}
	k \delta / 2 \ge \int_0^1 g_{\delta}(t) f_k(t) dt \ge \int_0^1 g_{\delta}(t) dX(t) = Y(\delta)/\delta,
\end{align*}
which implies that $k \delta^2 /2 \ge Y(\delta)$ for every $\delta >0$. But this contradicts Lemma~\ref{Lem:parabolictangent}, which says that a.s. $Y$ does not have a parabolic tangent at $t=0$. Consequently, $f_k$ is continuous at $t=0$. The same argument can also be applied to show the continuity of $f_k$ at $t=1$.
\hfill $\Box$
\end{prooftitle}
\medskip

\begin{lem}
\label{Lem:limknot}
Fix any $t \in (0,1)$. Denote by $\tau_k^-$ the right-most knot of $f_k$ on (0,t]. If such knot does not exist, then set $\tau_k^- = 0$. Similarly, denote by $\tau_k^+$ the left-most knot of $f_k$ on [t,1), and set $\tau_k^+=1$ if such knot does not exist. Then, for almost every $\omega \in \Omega$, there exists $K>0$ such that for every $k \ge K$, $\tau_k^- \neq 0$ and $\tau_k^+ \neq 1$. Here we suppressed the dependence of $K$ and $f_k$ (as well as $\tau_k^-$ and $\tau_k^+$) on $\omega$ (via $X$) in the notation.
\end{lem}
\begin{prooftitle}{of Lemma~\ref{Lem:limknot}}
Recall that $X(t)$ is said to be sample-bounded on $[0,1]$ if $\sup_{[0,1]}|X(t)|$ is finite for almost all $\omega \in \Omega$. Since both $X(t)$ and $Y(t)$ are sample bounded. Without loss of generality,
in the following, wee can fix $\omega \in \Omega$ and assume that $\sup_{[0,1]}|X(t)| < \infty$ and $\sup_{[0,1]}|Y(t)| < \infty$.

First, we show the existence of at least one knot on $(0,1)$. Note that the cubic polynomial $P_k$ with $P_k(0) = Y(0) = 0$, $P_k(1) = Y(1)$, $P_k^{(2)}(0) = P_k^{(2)}(1) = k$ can be expressed as
\begin{align*}
	P_k(s) = \frac{k}{2}s^2 + \left(Y(1) - \frac{k}{2}\right)s,
\end{align*}
Therefore, take for instance $s = 0.5$ and consider the event $P_k(0.5) \ge Y(0.5)$. This event can be reexpressed as
\[
\frac{1}{2}Y(1) - Y(0.5) \ge \frac{1}{8}k,
\]
which will eventually become false as $k \rightarrow \infty$. This is due to the fact that $Y(t)$ is sample bounded. In view of (i) and (v) of Lemma~\ref{Lem:limprop}, we conclude that $f_k$ has at least one knot in the open interval $(0,1)$ for sufficiently large $k$. 

Next, take $k$ large enough so that $f_k$ has one knot in $(0,1)$, which we denote by $\tau_k$. By (iii) of Lemma~\ref{Lem:limprop}, $H_k(\tau_k) = Y(\tau_k)$ and $H'_k(\tau_k) = X(\tau_k)$. Without loss of generality, 
we may assume that $\tau_k > t$. Now the cubic polynomial $P_k$ with $P_k(0) = Y(0) = 0$, $P_k(\tau_k) = Y(\tau_k)$, $P_k^{(2)}(0) = k$ and $P_k'(\tau_k) = X(\tau_k)$ can be expressed as 
\begin{align*}
	P_k(s) = a_{k1} s^3 + a_{k2} s^2 + a_{k3} s,
\end{align*}
where 
\begin{align*}
	a_{k1} &= \frac{X(\tau_k)}{2\tau_k^2} - \frac{Y(\tau_k)}{2\tau_k^3} - \frac{k}{4\tau_k};\\
	a_{k2} &= \frac{k}{2};\\
	a_{k3} &= \frac{3Y(\tau_k)}{2\tau_k} - \frac{X(\tau_k)}{2} - \frac{k\tau_k}{4}.
\end{align*}
By taking, say for example, $s = t/2$, it can then be verified that the event $P_k(t/2)\ge Y(t/2)$, which is equivalent to
\begin{align*}
	\frac{X(\tau_k)t^3}{16\tau_k^2} - \frac{Y(\tau_k)t^3}{16\tau_k^3} + \frac{3Y(\tau_k)t}{4\tau_k} - \frac{X(\tau_k)t}{4} \ge \frac{kt}{32\tau_k} (t - 2 \tau_k)^2,
\end{align*}
will eventually stop happening as $k \rightarrow \infty$. This is due to the sample boundedness of both $X(t)$ and $Y(t)$. Consequently, $\tau_k^- \neq 0$ for sufficiently large $k$. Furthermore, using essentially the same argument, one can also show that $\tau_k^+ \neq 1$ for large $k$, which completes the proof of this lemma.
\hfill $\Box$
\end{prooftitle}
\medskip

\begin{lem}
\label{Lem:limlowerbound}
For almost every $\omega \in \Omega$ (which determines $X(t)$ and $f_k$), we can find an $M > 0$ such that 
\[
\inf_{k \in \NN} \inf_{[0,1]} f_k > -M.
\]
\end{lem}
\begin{prooftitle}{of Lemma~\ref{Lem:limlowerbound}}
Fix any $\delta \in (0,1/2]$. In view of Lemma~\ref{Lem:limknot}, we may assume that there exist knots $\tau_k^-$ and $\tau_k^+$ on $(0,\delta]$ and $[1-\delta,1)$ respectively for sufficiently large $k$.  If $\inf_{[0,1]}f_k \ge 0$ for every $k > K$, then we are done. Otherwise, we focus on those $\inf_{[0,1]}f_k < 0$ and find $0< t_{k,1} < t_{k,2} < 1$ with $f_k(t_{k,1}) = f_k(t_{k,2}) = 0$. Note that the existence of $t_{k,1}$ and $t_{k,2}$ are guaranteed by Lemma~\ref{Lem:limprop} (v). In the following, we take $\delta = 1/12$ and consider two scenarios.
\begin{enumerate}[(a)]
\item $t_{k,2} - t_{k,1} < 2 \delta$. Let $a_{k,1} \in \partial f_k(t_{k,1})$ and $a_{k,2} \in \partial f_k(t_{k,2})$, where $\partial$ is the subgradient operator. Then Lemma~\ref{Lem:limprop}(v) implies that $a_{k,1} < 0$ and $a_{k,2} > 0$. Since both $a_{k,1}(s-t_{k,1})$ and $a_{k,2}(s-t_{k,2})$ can be regarded as supporting hyperplanes of $f_k$ due to its convexity, it follows that $f_k(s) \ge \max\big( a_{k,1}(s-t_{k,1}),a_{k,2}(s-t_{k,2}) \big)$. Set $C_k = - a_{k,1} a_{k,2} (t_{k,2} - t_{k,1})/(a_{k,2} - a_{k,1})$ (i.e. $-C_k$ is the value of the above hyperplanes at their intersection, which is negative), so $\inf_{[0,1]}f_k \ge -C_k$. Now
\begin{align*}
2 \sup_{t \in [0,1]}{|X(t)|} &\ge X(\tau_k^+) - X(\tau_k^-) = H_k'(\tau_k^+) - H_k'(\tau_k^-) = \int_{\tau_k^-}^{\tau_k^+} f_k(s) \, ds \\
	&\ge \int_{\tau_k^-}^{\tau_k^+} \max \big( a_{k,1}(s-t_{k,1}),a_{k,2}(s-t_{k,2}) \big) \, ds \\
	&\ge \int_{\delta}^{1-\delta} \max \big( a_{k,1}(s-t_{k,1}),a_{k,2}(s-t_{k,2}) \big) \, ds - 2\delta C_k \\
    &\ge \int_{\delta}^{1-\delta} \max \big( a_{k,1}(s-t_{k,1}),a_{k,2}(s-t_{k,2}), 0\big) \, ds \\
    & \qquad \qquad - \int_{\delta}^{1-\delta} \max \big( - a_{k,1}(s-t_{k,1}),- a_{k,2}(s-t_{k,2}), 0\big) ds - 2\delta C_k \\
    &\ge \inf_{u \in [0,1-4\delta]} \big\{|a_{k,1}| u^2 + |a_{k,2}|(1-4\delta - u)^2\big\}/2 - C_k (t_{k,2} - t_{k-1})/2  - 2 \delta C_k \\
	&\ge \frac{(1-4\delta)^2}{4\delta}C_k -\delta C_k - 2 \delta C_k =  \left(\frac{1}{4\delta} +\delta - 2\right) C_k \ge C_k.
\end{align*}
Consequently, a.s., $\sup_{k \in \NN} C_k < \infty$.
\item $t_{k,2} - t_{k,1} \ge 2 \delta$. Now consider $f_{k+} = \max(f_k,0)$. It follows that
\begin{align*}
	0 &\ge \phi(f_k) - \phi(f_{k+}) \\
	&=  \frac{1}{2} \int_0^1 \min(0,f_k(t))^2 \, dt - \int_0^1 \min(0,f_k(t))\, dX(t)
\end{align*}
Let $M_k = -\inf_{[0,1]} f_k$. By the convexity of $f_k$, we see that the first term is no smaller than $2 \delta M_k^2/3$. For the second term, we can again use Dudley's theorem and the fact that the following class
\begin{align*}
	\mathcal{G}_{\text{unit}} = \Big\{ g:[0,1]\rightarrow &[-1,0] \Big|  \exists \, 0 \le x_1 < x_2 \le 1 \\
	&\mbox{s.t. } g(s) \mbox{ is convex on } [x_1,x_2]; g(s) = 0, \forall s \in [0,x_1]\cup[x_2,1] \Big\}
\end{align*}
has entropy of order $\eta^{-1/2}$ in $L^2$ to argue that  
\begin{align*}
	 \sup_{g \in \mathcal{G}_{\text{unit}}} \left| \int_{[0,1]} g d X(t) \right| < \infty, \quad \mbox{a.s.}
\end{align*}
Therefore, the second term is at most $O(M_k)$.  Then we can use the argument in the proof of Lemma~\ref{Lem:limsolution} to establish that a.s. $\limsup_{k\rightarrow \infty} M_k < \infty$.
\end{enumerate}
\hfill $\Box$
\end{prooftitle}
\medskip

\begin{lem}
\label{Lem:limtight}
For any fixed $t \in (0,1)$ and almost every $\omega \in \Omega$,  $\sup_k |f_k(t)|$, $\sup_k |f_k^-(t)|$ and $\sup_k |f_k^+(t)|$ are bounded.
\end{lem}
\begin{prooftitle}{of Lemma~\ref{Lem:limtight}}
Let $\Delta = \min(t,1-t)/2$. In view of Lemma~\ref{Lem:limknot}, for sufficiently large $k$, we can assume that $f_k$ has at least one knot in $(0,t- \Delta]$, and one knot in $[t + \Delta,1)$. Denote these two points by $\tau_k^-$ and $\tau_k^+$ respectively. By the convexity of $f_k$, there exists some $c \in \RR$ such that $f_k(s) \ge c (s-t) + f_k(t)$ for every $s \in [0,1]$. It follows that
\begin{align*}
	2 \sup_{t \in [0,1]}{|X(t)|} &\ge X(\tau_k^+) - X(\tau_k^-) =  \int_{\tau_k^-}^{\tau_k^+} f_k(s) \, ds \\
	&\ge \int_{[\tau_k^-,t-\Delta] \cup [t+\Delta,\tau_k^+]} \Big(\inf_{[0,1]}f_k \Big) \, ds + \int_{(t-\Delta,t+\Delta)} \big\{c(s-t) + f_k(t) \big\} ds \ge -M_k + 2 \Delta f_k(t),
\end{align*}
where $M_k = - \inf_{[0,1]} f_k$. By Lemma~\ref{Lem:limlowerbound}, we see that for almost every $\omega \in \Omega$, $\limsup_{k \rightarrow \infty} f_k(t)$ is bounded. Combining this with the lower bound we established previously entails the boundedness of $\{f_k(t)\}_k$.

Next, note that both $\{f_k(t - \Delta)\}_k$ and $\{f_k(t + \Delta)\}_k$ are bounded. The boundedness of $\{f_k^-(t)\}_k$ and $\{f_k^+(t)\}_k$ immediately follows from the convexity of $f_k$ by utilizing
\begin{align*}
	\partial f_k(t) \subseteq \left[ \frac{f_k(t) - f_k(t-\Delta)}{\Delta}, \frac{f_k(t+\Delta) - f_k(t)}{\Delta} \right].
\end{align*}
\hfill $\Box$
\end{prooftitle}
\medskip

\begin{lem}
\label{Lem:limintegral}
For almost every $\omega \in \Omega$, both $\left\{\sup_{t \in [0,1]} | H_k(t)|\right\}_k$ and $\left\{\sup_{t \in [0,1]} | H_k'(t)|\right\}_k$ are bounded.
\end{lem}
\begin{prooftitle}{of Lemma~\ref{Lem:limintegral}}
By Lemma~\ref{Lem:limprop}(iv), for any $t \in [0,1]$,
\begin{align*}
	\min\Big(0, \inf_{k \in \NN} \inf_{[0,1]} f_k\Big) &\le t \inf_{[0,1]} f_k \le \int_0^t f_k(s) ds =  \int_0^1 f_k(s) ds - \int_t^1 f_k(s) ds\\
	&\le X(1) - X(0) - (1-t) \inf_{[0,1]} f_k \le X(1) - X(0) - \min\Big(0, \inf_{k \in \NN} \inf_{[0,1]} f_k\Big).
\end{align*}
Consequently, it follows from Lemma~\ref{Lem:limlowerbound} that $\left\{\sup_{t \in [0,1]}|\int_0^t f_k(s) ds|\right\}_k$ is bounded. Furthermore, Lemma~\ref{Lem:limknot} says that one can always find a knot $\tau_k^- \in (0,1)$ with $X(\tau_k^-) = H'_k(\tau_k^-)$ for all sufficiently large $k$. Thus, the boundedness of $\left\{\sup_{t \in [0,1]} | H_k'(t)|\right\}_k$ follows from the fact that
\begin{align*}
	\sup_{t \in [0,1]} | H_k'(t)| \le \sup_{t \in [0,1]} \left|X(\tau_k^-) - \int_0^{\tau_k^-} f_k(s)ds + \int_0^t f_k(s) ds\right| \le \sup_{t \in [0,1]}|X(t)| + 2\sup_{t \in [0,1]}\left|\int_0^t f_k(s) ds\right|.
\end{align*}
Finally, one can derive the sample boundedness of $\left\{\sup_{t \in [0,1]} | H_k(t)|\right\}_k$ by using the equality $H_k(t) = Y(0) + \int_0^t H_k'(s) ds$.
\hfill $\Box$
\end{prooftitle}
\medskip

\begin{lem}
\label{Lem:limequicontinuous}
For almost every $\omega \in \Omega$, both $\{H_k\}_k$ and $\{H_k'\}_k$ are uniformly equicontinuous on $[0,1]$. In fact, they are uniformly H\"{o}lder continuous with exponent less than $1/4$.
\end{lem}
\begin{prooftitle}{of Lemma~\ref{Lem:limequicontinuous}}
Here we only show that the family $\{H_k'\}_k$ is uniformly equicontinuous. Fix any $0 < \delta < 1$. For any $0 \le t_1 \le t_2 \le 1$ with $t_2 - t_1 < \delta$, by the convexity of $f_k$,
\begin{align*}
	H_k'(t_2) - H_k'(t_1) = \int_{t_1}^{t_2} f_k(s) ds \le \max\left(\int_{0}^{\delta} f_k(s) ds,\int_{1-\delta}^{1} f_k(s) ds \right).  
\end{align*}

In the following, we shall focus on $\int_{0}^{\delta} f_k(s) ds$. The term $\int_{1-\delta}^{1} f_k(s) ds$ can be handled in exactly the same fashion. By Lemma~\ref{Lem:limknot}, $\mathcal{S}(f_k)$ is non-empty for every $k > K$, where $K$ is a sufficiently large positive integer. Furthermore, in view of Lemma~\ref{Lem:limlowerbound}, we can assume that $\inf_{k \in \NN} \inf_{[0,1]} f_k \ge -M$ for some $M > 0$. 

Three scenarios are discussed in the following. Here we fix $\alpha \in (0,1/2)$.
\begin{enumerate}[(a)]
\item There exists at least one knot $\tau_k \in \mathcal{S}(f_k)$ with $\tau_k \in [\delta, \delta^\alpha]$. Then 
\begin{align*}
	\int_{0}^{\delta} f_k(s) ds &\le \int_{0}^{\tau_k} f_k(s) ds + M \delta^\alpha \\
	&\le X(\tau_k) - X(0) + M \delta^\alpha \le M' \delta^{{\alpha}^2} + M \delta^\alpha,
\end{align*}
for some $M' > 0$ (which is the $\alpha$-H\"{o}lder constant of this particular realization of $X(t)$). The last line follows from Lemma~\ref{Lem:limprop}(iv) and the fact that $X(t)$ is $\alpha$-H\"{o}lder-continuous, so 
\[
|X(\tau_k)-X(0)| \le M' |\tau_k|^\alpha \le M'(\delta^\alpha)^\alpha = M'\delta^{{\alpha}^2}.
\] 

\item $\mathcal{S}(f_k) \cap (0, \delta^\alpha] = \emptyset$. Let $\tau_k \in \mathcal{S}(f_k)$ be the left-most knot in $[\delta^\alpha,1)$. Then $f_k$ is linear on $[0,\tau_k]$ and $f_k(\tau_k) \ge -M$. Note that
\begin{align*}
	M' {\tau_k}^\alpha \ge X(\tau_k) - X(0) \ge \int_0^{\tau_k} f_k(s) ds = \frac{f_k(0) + f_k(\tau_k)}{2} \tau_k  \ge (f_k(0) - M) \tau_k /2, 
\end{align*}
so $f_k(0) \le 2M' \tau_k^{\alpha-1} + M \le 2M' \delta^{\alpha^2 -\alpha} + M$. It now follows that 
\begin{align*}
	\int_{0}^{\delta} f_k(s) ds \le (2M' \delta^{\alpha^2 -\alpha} + M )\delta \le (2M'+M) \delta^\alpha.
\end{align*}
We remark that since $f_k(0)=k$, the above conclusion also implies that $\mathcal{S}(f_k) \cap (0, \delta^\alpha] \neq \emptyset$ for all sufficiently large $k$.

\item $\mathcal{S}(f_k) \cap [\delta, \delta^\alpha) = \emptyset$, but $\mathcal{S}(f_k) \cap (0,\delta] \neq \emptyset$. Let $\tau_k^-$ be the right-most knot in $(0,\delta]$ and $\tau_k^+$ be the left-most knot in $[\delta^\alpha,1)$. As a convention, we set $\tau_k^+ = 1$ if such a knot does not exist. Note that $f_k$ is linear on $[\tau_k^-,\tau_k^+]$, so we can use essentially the same argument as above to see that $\int_{\tau_k^-}^\delta f_k(s) ds \le (2M'+M) \delta^\alpha$. Finally,
\begin{align*}
	\int_{0}^{\delta} f_k(s) ds \le X(\tau_k^-)-X(0) + \int_{\tau_k^-}^{\delta} f_k(s) ds \le M' \delta^\alpha + (2M'+M) \delta^\alpha = (3M'+M)\delta^\alpha.
\end{align*}
\end{enumerate}

To finish the proof, we shall apply Lemma~\ref{Lem:limlowerbound} to verify that $H_k'(t_2) - H_k'(t_1) \ge -M \delta$.

In addition, combining results from all the scenarios listed above, we see that $\{H_k'\}_k$ is uniformly H\"{o}lder continuous with exponent $\alpha^2$ ($< 1/4$).
\hfill $\Box$
\end{prooftitle}
\medskip

\begin{prooftitle}{of Theorem~\ref{Thm:triangular_limit}}
For every $m \in \NN$ with $m \ge 3$, define the following norms
\begin{align*}
	\lVert H \rVert_m = \sup_{t \in [0, 1]} |H(t)| + \sup_{t \in [0, 1]} |H'(t)| + \sup_{t \in [1/m, 1-1/m]}|H^{(2)}(t)|.
\end{align*}

First, we show the existence of such a function for almost all $\omega \in \Omega$ by construction. Fix $\omega$ (thus we focus on a particular realization of $X(t)$ but suppress its dependence on $\omega$ in the notation). Let $H_k$ be the function satisfying $H_k^{(2)} = f_k$, $H_k(0) = 0$ and $H_k(1) = Y(1)$. We claim that the sequence $H_k$ admits a convergent subsequence in the topology induced by the norm $\lVert \cdot \rVert_m$. 

By Lemma~\ref{Lem:limtight}, we may assume that for $t = 1/m$ and $t= 1-1/m$, $\{f_k(t)\}_k$, $\{f_k^-(t)\}_k$ and $\{f_k^+(t)\}_k$ are bounded. By the convexity, $\{f_k\}_k$ have uniformly bounded derivatives on $[1/m,1-1/m]$ so are uniformly bounded and equicontinuous. Therefore, the Arzel\`{a}--Ascoli theorem guarantees that the sequence $f_k$ has a convergent subsequence $f_{k_l}$ in the supremum metric on $[1/m, 1-1/m]$. Extracting further subsequences if necessary, one can get $f_{k_l}$ converging in the topology induced by the $L^\infty$ norm on $[1/m,1-1/m]$ for $m = 3,4,\ldots$. 

Now by Lemma~\ref{Lem:limintegral} and Lemma~\ref{Lem:limequicontinuous}, we can assume that $\{H_k\}_k$ and $\{H_k'\}_k$ are bounded and equicontinuous on $[0,1]$. By the Arzel\`{a}--Ascoli theorem again, we are able to extract further subsequences if necessary to make $H_{k_l}$ converge in the topology induced by the norms $\lVert \cdot \rVert_m$ for $m = 3,4,\ldots$. We denote the function that $H_{k_l}$ converges to by $H$. 

In the following, we show that $H$ has the properties listed in the statement of the theorem. 
\begin{enumerate}
\item[(1)] Because $H_{k}(t) \ge Y(t)$ for every $k \in \NN$ and every $t \in [0,1]$, $H(t) \ge Y(t)$ for $t \in [0,1]$.
\item[(2)] $H^{(2)}$ is convex on $(0,1)$ since every $f_k$ is convex.
\item[(3) -- (4)] Since $H_k(0) = Y(0)$ and $H_k(1) = Y(1)$ for all $k \in \NN$, we have $H(0) = Y(0)$ and $H(1) = Y(1)$. Now let $\tau_k^- = \inf \{\tau_k: \tau \in \mathcal{S}(f_k)\}$ and $\tau_k^+ = \sup \{\tau_k: \tau \in \mathcal{S}(f_k)\}$. In light of Lemma~\ref{Lem:limknot}, one can assume that $\lim_{k \rightarrow \infty} \tau_k^- = 0$ and $\lim_{k \rightarrow \infty} \tau_k^+ = 1$. The sample continuity of $X(t)$, together with the property of knots (see (iii) of Lemma~\ref{Lem:limprop}), entails that $H'_k(\tau_k^-) \rightarrow X(0)$ and $H'_k(\tau_k^+) \rightarrow X(1)$ as $k \rightarrow \infty$. Finally, we use the uniform convergence and the continuity of $H'$ to establish $H'(0) = X(0)$ and $H'(1) = X(1)$.
\item[(5)] It follows from Lemma~\ref{Lem:limprop}(ii) that $\int_{[1/m,1-1/m]} \big(H(t)-Y(t)\big)\, dH^{(3)}(t) = 0$. Now let $m \rightarrow \infty$ to see the required property.
\end{enumerate}
This completes the proof of existence. 

It remains to show the uniqueness of $H$. Suppose that there are $H_1$ and $H_2$ satisfying Conditions (1) -- (5) listed in the statement of Theorem~\ref{Thm:triangular_limit}. For notational convenience, we write $h_1 = H_1^{(2)}$ and $h_2 = H_2^{(2)}$. Then,
\begin{align*}
& \bigg\{\frac{1}{2}\int_0^1 h_1^2(t) dt - \int_0^1 h_1(t)dX(t)\bigg\} -  \bigg\{\frac{1}{2} \int_0^1 h_2^2(t) dt - \int_0^1 h_2(t)dX(t)\bigg\} \\
& = \frac{1}{2} \int_0^1 \Big(h_1(t)-h_2(t)+h_2(t)\Big)^2 dt - \int_0^1 \Big(h_1(t)-h_2(t)+h_2(t)\Big)dX(t)   \\
& \qquad -\frac{1}{2} \int_0^1 h_2^2(t) dt + \int_0^1 h_2(t)dX(t) \\
& = \frac{1}{2} \int_0^1 \big(h_1(t) - h_2(t)\big)^2 dt + \bigg\{ \int_0^1 \big(h_1(t)-h_2(t)\big)h_2(t)dt - \int_0^1 \big(h_1(t)-h_2(t)\big)dX(t)\bigg\} \\
& \ge \frac{1}{2} \int_0^1 \big(h_1(t) - h_2(t)\big)^2 dt,
\end{align*}
where we used Conditions (1) -- (5) of $H_2$ to derive the last inequality.
By swapping $H_1$ and $H_2$, we further obtain the following inequality
\[
\bigg\{\frac{1}{2} \int_0^1 h_2^2(t) dt - \int_0^1 h_2(t)dX(t)\bigg\} -  \bigg\{\frac{1}{2} \int_0^1 h_1^2(t) dt - \int_0^1 h_1(t)dX(t)\bigg\} \ge \frac{1}{2} \int_0^1 \big(h_2(t) - h_1(t)\big)^2 dt.
\] 
Adding together the above two inequalities yields $0 \ge \int_0^1 \big(h_1(t) - h_2(t)\big)^2 dt$,
which implies the uniqueness of $H^{(2)}$ on $(0,1)$. The uniqueness of $H$ then follows from its third condition.

The proof for the second part (i.e. the existence and uniqueness of $\tilde{H}$) is similar and is therefore omitted.
\hfill $\Box$
\end{prooftitle}
\medskip

\begin{prooftitle}{of Corollary~\ref{Cor:triangular_limit}}
We can easily verify the existence of such a function by using the same construction in the proof of Theorem~\ref{Thm:triangular_limit}. In particular, if $Y(t)$ does not have parabolic tangents at both $t = 0$ and $t = 1$ (which happens a.s. according to Lemma~\ref{Lem:parabolictangent}), then $f_k(0) \rightarrow \infty$ and $f_k(1) \rightarrow \infty$ as $k \rightarrow \infty$ .

On the other hand, if $H^{(2)}(0^+) \rightarrow \infty$, there must be a sequence of knots $\tau_1,\tau_2,\ldots$ of $H^{(2)}$ with $\lim_{j \rightarrow \infty} \tau_j = 0$. In views of Conditions~(1), (2) and (5), one necessarily has $H(\tau_j) = Y(\tau_j)$ and $H'(\tau_j) = X(\tau_j)$ for every $j$. The fact that $H$, $H'$, $Y$ and $X$ are all continuous entails that $H(0) = Y(0)$ and $H'(0) = X(0)$. Consequently, Condition~(3') implies Condition~(3). We now apply the same argument to $H^{(2)}(1^-)$ to conclude that Condition~(4') implies Condition~(4). Hence, in view of Theorem~\ref{Thm:triangular_limit}, $H$ is unique. 
\hfill $\Box$
\end{prooftitle}
\medskip

\subsection{Appendix II: pointwise adaptation for cases (B) and (C)}
The following three lemmas are required to prove Theorem~\ref{Thm:LS_adapt_2}.
\begin{lem}
\label{Lem:LS_adapt_precalculation}
For any $\alpha>1$, $\inf_{k \in [0,1]} \Big\{ 4k^{\alpha+2}+(1+k)^{\alpha+1}(\alpha-2k)\Big\} > 0$.
\end{lem}
\begin{prooftitle}{of Lemma~\ref{Lem:LS_adapt_precalculation}}
It suffices to show that for any $k \in [0,1]$, 
\[
4k^{\alpha+2}+(1+k)^{\alpha+1}(\alpha-2k) > 0.
\]

First, it is easy to check that the above inequality holds true when $k=0$. In the case of $k > 0$, we can restate the inequality to be proved as 
\begin{align}
\label{Eq:LS_adapt_precalculation1}
4 \Big(\frac{k}{1+k}\Big)^{\alpha+1} > 2 - \frac{\alpha}{k}.
\end{align}

Next, we define $m = k/(1+k) \in [0,1/2]$, so that (\ref{Eq:LS_adapt_precalculation1}) can be rewritten as
\begin{align}
\label{Eq:LS_adapt_precalculation2}
4 m^{\alpha+1} + \frac{\alpha}{m} > 2+\alpha.
\end{align}
Note that 
\[
\inf_{m \in [0,1/2]} \Big( 4 m^{\alpha+1} + \frac{\alpha}{m} \Big) \ge \frac{\alpha(\alpha+2)}{(\alpha+1)} \Big(\frac{4(\alpha+1)}{\alpha}\Big)^{1/(\alpha+2)} =  (\alpha+2) \bigg\{4\Big(1 - \frac{1}{\alpha+1}\Big)^{\alpha+1}\bigg\}^{1/(\alpha+2)}.
\]
The inequality (\ref{Eq:LS_adapt_precalculation2}) now follows easily from that fact that $\Big(1 - \frac{1}{\alpha+1}\Big)^{\alpha+1} > 1/4$ for any $\alpha>1$.
\hfill $\Box$
\end{prooftitle}
\medskip

\begin{lem}
\label{Lem:LS_adapt_calculation}
Let $f_0(t) = f_0(x_0) K_1 (t-x_0) + K_2(t-x_0)^\alpha \mathbf{1}_{[x_0,\infty)}(t)$ for every $t \in [x_0 - \delta, x_0 + \delta]$ for some $\delta \in (0,x_0)$ small, with $K_1 < 0$, $K_2 > 0$ and $\alpha > 1$. Then for any $x_0 - \delta \le \tau^- \le x_0 \le \tau^+ \le x_0 + \delta$,
\[
\int_{\tau^-}^{\tau} \bigg\{\frac{1}{2}(\tau^-+\tau) - t \bigg\} f_0(t)dt 
+ \int_{\tau}^{\tau^+} \bigg\{t - \frac{1}{2}(\tau+\tau^+) \bigg\} f_0(t)dt \ge K_{f_0} (\tau^+ -x_0)^{\alpha+1} (\tau^+ - \tau^-),
\]
where $\tau = (\tau^- + \tau^+)/2$, and where $K_{f_0} > 0$ is a constant that only depends on $f_0$. 
\end{lem}
\begin{prooftitle}{of Lemma~\ref{Lem:LS_adapt_calculation}}
First, it is easy to check that
\begin{align}
\notag &\int_{\tau^-}^{\tau} \bigg\{\frac{1}{2}(\tau^-+\tau) - t \bigg\} f_0(t)dt 
+ \int_{\tau}^{\tau^+} \bigg\{t - \frac{1}{2}(\tau+\tau^+) \bigg\} f_0(t)dt \\
\label{Eq:LS_adapt_calculation1} 
&= K_2 \int_{\tau^-}^{\tau} \bigg\{\frac{1}{2}(\tau^-+\tau) - t \bigg\} (t-x_0)^\alpha  \mathbf{1}_{[x_0,\infty)}(t) dt 
+ K_2\int_{\tau}^{\tau^+} \bigg\{t - \frac{1}{2}(\tau+\tau^+) \bigg\} (t-x_0)^\alpha  \mathbf{1}_{[x_0,\infty)}(t) dt.
\end{align}

If $\tau \le x_0$, then (\ref{Eq:LS_adapt_calculation1}) can be expressed as
\begin{align*}
K_2 \int_{0}^{\tau^+ - x_0} \bigg\{t - \frac{1}{2}(\tau+\tau^+) + x_0 \bigg\} t^\alpha dt &= \frac{K_2(\tau^+ -x_0)^{\alpha+1}\{(\tau^+ - \tau)\alpha +  2(x_0 - \tau)\}}{2(\alpha+1)(\alpha+2)} \\
& \ge \frac{K_2\alpha}{4(\alpha+1)(\alpha+2)}(\tau^+ -x_0)^{\alpha+1}(\tau^+ - \tau^-).
\end{align*}

On the other hand, if $\tau > x_0$, then after some elementary calculations, we can show that (\ref{Eq:LS_adapt_calculation1}) is equal to
\begin{align}
\label{Eq:LS_adapt_calculation2}
\frac{K_2 [4(\tau-x_0)^{\alpha+2} + (\tau^+ - x_0)^{\alpha+1} \{ \alpha(\tau^+ - \tau) - 2(\tau-x_0) \}]}{2(\alpha+1)(\alpha+2)}.
\end{align}
Denote by $k = (\tau-x_0)/(\tau^+ - \tau)$, so that (\ref{Eq:LS_adapt_calculation2}) can be rewritten as 
\[
\frac{K_2}{{2(\alpha+1)(\alpha+2)}} \big\{4k^{\alpha+2}+(1+k)^{\alpha+1}(\alpha-2k)\big\} (\tau^+ - \tau)^{\alpha+2} \ge C_{f_0}(\tau^+ - \tau)^{\alpha+2},
\]
where $C_{f_0} > 0$ is a constant that only depends on $f_0$, and where we applied Lemma~\ref{Lem:LS_adapt_precalculation} with the fact that $k \in [0,1]$ to derive the above displayed equation. 
Consequently, by setting $K_{f_0} = C_{f_0}/2^{\alpha+2}$, it is straightforward to check that (\ref{Eq:LS_adapt_calculation2}) is greater than or equal to $K_{f_0}(\tau^+ - x_0)^{\alpha+1}(\tau^+ - \tau^-)$.
\hfill $\Box$
\end{prooftitle}
\medskip

\begin{lem}
\label{Lem:LS_adapt_VC}
Let $\mathcal{F}$ be a collection of functions defined on $[x_0 - \delta, x_0 + \delta]$, with $\delta > 0$ 
small. Suppose that for a fixed $x \in [x_0 - \delta, x_0 + \delta]$ and 
every $0 < R \le R_0$ such that $[x, x+ R] \subseteq [x_0 - \delta, x_0 + \delta]$, the collection
\begin{align*}
\mathcal{F}_{x,R} 
= \left \{f_{x,y} = f 1_{[x,y]}, \ \ \ f \in \mathcal{F}, \ \ x \le y \le x+R \right \}
\end{align*}
admits an envelope $F_{x,R}$ such that 
\[
\mathbb{E}F^2_{x,R}(X_1) \leq K R^{2d-1}
\]
for some $d \ge 1/2$ fixed and $K>0$ depending only on $x_0$ and $\delta$, where $X_1 \sim F_0$.
Moreover, suppose that
\[
\sup_ Q \int_0^1 \sqrt{\log N(\eta \Vert F_{x,R} \Vert_{Q,2}, 
\mathcal{F}_{x,R}, L_2(Q))}d\eta < \infty. 
\]
Then, for every $\epsilon > 0$ and $s_0 > 0$, there exist a random variable $M_n$ 
of $O_p (1)$ such that
\begin{align*}
  \sup_{s \ge s_0} \; \sup_{x \le y \le x+R_0} n^{(s+d)/(2s+1)} \max \bigg(\sup_{f_{x,y} \in \mathcal{F}_{x,y-x}} \bigg| \int f_{x, y} d(\mathbb{F}_n - F_0)\bigg| - \epsilon |y-x|^{s+d}, 0\bigg)  = M_n.
\end{align*} 
\end{lem}
\begin{prooftitle}{of Lemma~\ref{Lem:LS_adapt_VC}}
This lemma slightly generalizes Lemma~A.1 of \citet{BalabdaouiWellner2007}. Its proof proceeds as in \citet{BalabdaouiWellner2007} with minor modifications, so is omitted for brevity. We remark that here only the collection of functions defined on $[x,x+R]$ are considered. For functions on $[x-R,x]$, an analogous version of this lemma also holds true by symmetry.
\hfill $\Box$
\end{prooftitle}
\medskip

\begin{prooftitle}{of Theorem~\ref{Thm:LS_adapt_2}}
Here we only consider case (B). Case (C) can be handled similarly by symmetry.

Suppose that (B) holds true for every $t \in [x_0 - 2\delta, x_0 + 2\delta]$ for some fixed $\delta > 0$. Let $\tau_n^-  = \max\{t \in \mathcal{S}(\hat{f}_n), t < x_0\}$ and $\tau_n^+ = \min\{t \in \mathcal{S}(\hat{f}_n), t \ge x_0\}$. Note that consistency of $\hat{f}_n$ implies that $\tau_n^+ - x_0 = o_p(1)$. Moreover, since $f_0$ is linear on $[x_0 - \delta, x_0]$, we can proceed as in the proofs of Theorem~\ref{Thm:LS_triangular_rate} and Corollary~\ref{Cor:LS_triangular_uniform_rate_derivative} to show that
\begin{align}
\label{Eq:LS_adapt_2_proof_0}
\inf_{t \in [x_0 - \delta, x_0]} \min\big(\hat{f}_n(t) - f_0(t),0\big)  = O_p(n^{-1/2})
\end{align}
and 
\[
\sup_{t \in [x_0 - \delta, x_0]} |\hat{f}_n^-(t) - f_0'(t)|  = O_p(n^{-1/2}),
\]
where $\hat{f}_n^-$ is the left derivative of $\hat{f}_n$. Therefore, it suffices to only consider the behavior of $\max\big(\hat{f}_n(x_0) - f_0(x_0),0\big)$. 

The proof can be divided into four parts.
\begin{enumerate}[(i)]
\item Suppose that $x_0 - \tau_n^- \ge n^{-1/(4\alpha+2)}$. Denote by $\tau_n^* = \max(\tau_n^-, x_0-\delta)$. Then because both $\hat{f}_n$ and $f_0$ are linear on $[\tau_n^*,x_0]$, 
\[
\big(\hat{F}_n(x_0) - F_0(x_0)\big) - \big(\hat{F}_n(\tau_n^*) - F_0(\tau_n^*)\big) = (x_0 - \tau_n^-)\big(\hat{f}_n(x_0) - f_0(x_0) + \hat{f}_n(\tau_n^*) - f_0(\tau_n^*)\big)/2.
\]
Marshall's lemma entails that the left-hand side of the above equality is $O_p(n^{-1/2})$. Furthermore, it has been shown that 
\[
\min\big(\hat{f}_n(\tau_n^*) - f_0(\tau_n^*),0) = O_p(n^{-1/2}).
\]
Rearranging the above equation yields 
\[
\max\big(\hat{f}_n(x_0) - f_0(x_0),0\big) = O_p(n^{-\alpha/(2\alpha+1)}).
\]

\item Suppose that $\tau_n^+ - x_0 < n^{-1/(2\alpha+1)}$. First, we modify Lemma~4.2 of \citet{GJW2001b} to the following: 
\begin{enumerate}
\item[] Let $\xi_n$ be an arbitrary sequence of numbers converging to $x_0$. Define $\eta_n^- = \max\{t \in \mathcal{S}(\hat{f}_n): t < \xi_n\}$ and $\eta_n^+ = \min\{t \in \mathcal{S}(\hat{f}_n): t \ge \xi_n\}$. Then, $\max(\eta_n^+,\tau^+) - \max(\eta_n^-, \tau_n^+) = O_p(n^{-\alpha/(2\alpha+1)})$. 
\end{enumerate}
Here one can apply Lemma~A.1 of \citet{BalabdaouiWellner2007} with $k = \alpha$ and $d = 2$ to verify the above extension. By proceeding as in the proof of Lemma~4.3 of \citet{GJW2001b}, one can verify that there exists some $\tau_n^{+++} > \tau_n^+$ such that $\tau_n^{+++} - x_0 = O_p(n^{-1/(2\alpha+1)})$ and 
\[
\inf_{t \in [\tau_n^{+}, \tau_n^{+++}]} |\hat{f}_n(t) - f_0(t)| = O_p(n^{-\alpha/(2\alpha+1)}).
\]
Since $|\hat{f}_n^-(x_0)-f_0'(x_0)| = O_p(n^{-1/2})$, it follows that
\begin{align*}
&\max(\hat{f}_n(x_0) - f_0(x_0),0) \\
&\le \inf_{t \in [\tau_n^{+}, \tau_n^{+++}]} |\hat{f}_n(t) - f_0(t)| - \hat{f}_n^-(x_0) (\tau_n^{+++}-x_0) - f_0(\tau_n^{+++}) + f_0(x_0) \\
&\le O_p(n^{-\alpha/(2\alpha+1)}) - (\hat{f}_n^-(x_0)-f_0'(x_0)) (\tau_n^{+++}-x_0) + K_2(\tau_n^{+++}-x_0)^\alpha = O_p(n^{-\alpha/(2\alpha+1)}).
\end{align*}
In fact, it is easy to see that the above conclusion still holds if we change our assumption from $\tau_n^+ - x_0 < n^{-1/(2\alpha+1)}$ to $\tau_n^+ - x_0 = O_p(n^{-1/(2\alpha+1)})$.

\item Suppose that $x_0 - \tau_n^- \le \tau_n^+ - x_0$. In view of the discussion in (ii), it suffices to show that $\tau_n^+ - x_0 = O_p(n^{-1/(2\alpha+1)})$. 
Consider the behavior of $\int_0^{t} \hat{F}_n(s) ds$ (as a function of $t$) at the middle point $\tau_n = (\tau_n^- + \tau_n^+)/2$. Let
\begin{align*}
E_{1n} &= \int_{\tau_n^-}^{\tau_n} \bigg\{\frac{1}{2}(\tau_n^-+\tau_n) - t \bigg\} d\big (F_0(t) - \mathbb{F}_n(t)\big) 
+ \int_{\tau_n}^{\tau_n^+} \bigg\{t - \frac{1}{2}(\tau_n+\tau_n^+) \bigg\} d\big(F_0(t) - \mathbb{F}_n(t)\big); \\
E_{2n} &= \int_{\tau_n^-}^{\tau_n} \bigg\{\frac{1}{2}(\tau_n^-+\tau_n) - t \bigg\} f_0(t)dt 
+ \int_{\tau_n}^{\tau_n^+} \bigg\{t - \frac{1}{2}(\tau_n+\tau_n^+) \bigg\} f_0(t)dt. 
\end{align*}
It was shown in Lemma~4.2 of \citet{GJW2001b} that 
\begin{align}
\label{Eq:LS_adapt_2_proof_1}
E_{1n} \ge E_{2n}.
\end{align}
In view of Lemma~A.1 of \citet{BalabdaouiWellner2007} with $k = \alpha$ and $d = 2$, 
\begin{align}
\label{Eq:LS_adapt_2_proof_2}
E_{1n} \le \epsilon (\tau_n^+ - \tau_n^-)^{\alpha+2} + O_p(n^{-(\alpha+2)/(2\alpha+1)}) \le \epsilon \{2(\tau_n^+ - x_0)\}^{\alpha+2} + O_p(n^{-(\alpha+2)/(2\alpha+1)})
\end{align}
for any $\epsilon > 0$. On the other hand, by Lemma~\ref{Lem:LS_adapt_calculation}, we have 
\begin{align}
\label{Eq:LS_adapt_2_proof_3}
E_{2n} \ge K_{f_0} (\tau_n^+ - x_0)^{\alpha+2}
\end{align}
for some constant $K_{f_0} > 0$ that only depends on the density function $f_0$. Combining (\ref{Eq:LS_adapt_2_proof_1}), (\ref{Eq:LS_adapt_2_proof_2}) and (\ref{Eq:LS_adapt_2_proof_3}) together yields $\tau_n^+ - x_0 = O_p(n^{-1/(2\alpha+1)})$.

\item Finally, suppose that $x_0 - \tau_n^- < n^{-1/(4\alpha+2)}$, $\tau_n^+ - x_0 \ge n^{-1/(2\alpha+1)}$ and $x_0 - \tau_n^- > \tau_n^+ - x_0$. Our first aim here is to show that
\begin{align}
\label{Eq:LS_adapt_2_proof_4}
\frac{(\tau_n^+ - x_0)^{\alpha+1}}{\tau_n^+ - \tau_n^-} = O_p(n^{-\alpha/(2\alpha+1)}).
\end{align}
Let 
\begin{align}
\label{Eq:LS_adapt_2_proof_5}
s_n =(\alpha+1) \frac{\log (\tau_n^+ - x_0)}{ \log (\tau_n^+ - \tau_n^-)} - 1.
\end{align}
Since $\tau_n^+ \stackrel{a.s.}{\rightarrow} x_0$ and $x_0 - \tau_n^- < n^{-1/(4\alpha+2)}$, $0 < \tau_n^+ - x_0 < \tau_n^+ - \tau_n^- \le 1$ for large $n$. This implies that $s_n \ge \alpha$ for sufficiently large $n$. 
It now follows from Lemma~\ref{Lem:LS_adapt_VC} (with $d = 2$ and $s_0 = \alpha$) that
\[
n^{(s_n+2)/(2s_n+1)} \max \Big(E_{1n} - \epsilon (\tau_n^+ - \tau_n^-)^{s_n+2}, 0\Big)  = O_p(1). 
\]
for any small $\epsilon > 0$.  On the other hand, in view of Lemma~\ref{Lem:LS_adapt_calculation} and identity (\ref{Eq:LS_adapt_2_proof_5}), some elementary calculations yield
\[
E_{2n} \ge K_{f_0}(\tau_n^+ - \tau_n^-)^{s_n+2}.
\]
Plugging the above two equations into (\ref{Eq:LS_adapt_2_proof_1}) entails that $\tau_n^+ - \tau_n^- \le n^{-1/(2s_n+1)} O_p(1)$.
As a result, 
\[
\frac{(\tau_n^+ - x_0)^{\alpha+1}}{\tau_n^+ - \tau_n^-} = (\tau_n^+ - \tau_n^-)^{s_n} = n^{-s_n/(2s_n+1)} O_p(1) \le O_p(n^{-\alpha/(2\alpha+1)}).
\]
The rest of the proof is similar to (b) in the proof of Theorem~\ref{Thm:LS_adapt_1}. By the law of the iterated logarithm for local empirical processes,
\begin{align}
\notag&\sqrt{\tau_n^+ - \tau_n^-} O_p\big(n^{-1/2} \sqrt{\log\log n}\big) \\
\notag&=\Big(\mathbb{F}_n(\tau_n^+) - F_0(\tau_n^+)\Big) - \Big(\mathbb{F}_n(\tau_n^-) - F_0(\tau_n^-)\Big) \\
\notag&= \int_{\tau_n^-}^{\tau_n^+} \big(\hat{f}_n(t) - f_0(t)\big)dt \\
\notag&= \left\{ \frac{\tau_n^+ - \tau_n^-}{2(x_0 - \tau_n^-)}(\hat{f}_n(x_0) - f_0(x_0)) + \Big(1 - \frac{\tau_n^+ - \tau_n^-}{2(x_0 - \tau_n^-)}\Big) (\hat{f}_n(\tau_n^-) - f_0(\tau_n^-))\right\}(\tau_n^+ - \tau_n^-) \\
\notag & \quad - \frac{K_2}{\alpha+1} (\tau_n^+ - x_0)^{\alpha+1},
\end{align}
where the last equality follows from the linearity of $\hat{f}_n(t)$ and $f_0(t) - K_2 (t-x_0)^\alpha  \mathbf{1}_{[x_0,\infty)}(t)$ on $[\tau_n^-,\tau_n^+]$.
Since our assumption in (iv) guarantees that $1 \le (\tau_n^+ - \tau_n^-)/(x_0 - \tau_n^-) \le 2$, rearranging the terms in the above displayed equation leads to
\begin{align}
\notag&\max\Big(\hat{f}_n(x_0) - f_0(x_0), 0 \Big) \\
\notag&\le 2\bigg\{ (\tau_n^+ - \tau_n^-)^{-1/2} O_p\big(n^{-1/2}\sqrt{\log\log n}\big) \\
\label{Eq:LS_adapt_2_proof_6}
 & \qquad\qquad-  \frac{1}{2}\min(\hat{f}_n(\tau_n^-) - f_0(\tau_n^-),0) + \frac{K_2}{\alpha+1} \frac{(\tau_n^+ - x_0)^{\alpha+1}}{\tau_n^+ - \tau_n^-}\bigg\}.
\end{align}
Finally, as $\tau_n^+ - \tau_n^- \ge n^{-1/(2\alpha+1)}$, we can plug (\ref{Eq:LS_adapt_2_proof_0}) and (\ref{Eq:LS_adapt_2_proof_4}) into (\ref{Eq:LS_adapt_2_proof_6}) to verify 
\[
\max\Big(\hat{f}_n(x_0) - f_0(x_0), 0 \Big) \le O_p\big(n^{-\alpha/(2\alpha+1)}\sqrt{\log\log n}\big).
\]
\end{enumerate}
\hfill $\Box$
\end{prooftitle}
\medskip

\section*{Acknowledgments}
The first author is grateful to Richard Samworth for helpful conversations. 
The second author owes thanks to Tony Cai and Mark Low for
questions concerning the problems addressed in Section 2.3. He also
owes thanks to Fadoua Balabdaoui and Hanna Jankowski for conversations 
and initial work on the problems addressed in Sections 2.1 and
2.2. Part of this work was completed while the first author visited the
Statistics Department at the University of Washington. The first author was 
supported by EPSRC grant EP/J017213/1. The second author was 
supported in part by NSF grant DMS-1104832 and by NI-AID grant
2R01 AI291968-04.



\begin{thebibliography}{75}
\bibitem[{Balabdaoui(2007)}]{Balabdaoui2007}
Balabdaoui, F. (2007) Consistent estimation of a convex density at the origin.
\newblock \emph{Math. Methods Statist.}, \textbf{16}, 77--95. 

\bibitem[{Balabdaoui(2014)}]{Balabdaoui2014}
Balabdaoui, F. (2014) Global convergence of the log-concave MLE when the truth is geometric. 
\newblock \emph{J. Nonparametric. Statist.}, \textbf{26}, 21--59. 

\bibitem[{Balabdaoui\emph{ et.al.}(2013)}]{BJRP2013}
Balabdaoui, F., Jankowski, H., Rufibach, K. and Pavlides, M.(2013).
Asymptotic distribution of the discrete log-concave MLE and related applications.
\newblock \emph{J. Roy. Statist. Soc. Ser. B}, \textbf{75}, 769--790. 

\bibitem[{Balabdaoui and Rufibach(2008)}]{BalabdaouiRufibach2008}
Balabdaoui, F. and Rufibach, K. (2008) A second Marshall inequality in convex estimation
\newblock \emph{Statist. Probab. Lett.}, \textbf{78}, 118--126.

\bibitem[{Balabdaoui, Rufibach and Wellner(2009)}]{BRW2009}
Balabdaoui, F., Rufibach, K. and Wellner, J.~A. (2009) Limit distribution
  theory for maximum likelihood estimation of a log-concave density.
\newblock \emph{Ann. Statist.}, \textbf{37}, 1299--1331.

\bibitem[{Balabdaoui and Wellner(2007)}]{BalabdaouiWellner2007}
Balabdaoui, F. and Wellner, J.~A. (2007) Estimation of a $k$-monotone 
density: limit distribution theory and the spline connection.
\newblock \emph{Ann. Statist.}, \textbf{35}, 2536--2564.

\bibitem[{Balabdaoui and Wellner(2014)}]{BalabdaouiWellner2014}
Balabdaoui, F. and Wellner, J.~A. (2014) Chernoff's density is log-concave.
\newblock \emph{Bernoulli}, \textbf{20}, 231--244.

\bibitem[{Cai and Low(2015)}]{CaiLow2015}
Cai, T. T. and Low, M. (2015) 
A framework for estimation of convex functions.
\newblock \emph{Statist. Sinica.}, \textbf{25}, 423-456.

\bibitem[{Cator(2011)}]{Cator2011}
Cator, E. (2011) 
Adaptivity and optimality of the monotone least-squares estimator
\newblock \emph{Bernoulli}, \textbf{17}, 714-735.

\bibitem[{Carolan and Dykstra(1999)}]{CarolanDykstra1999}
Carolan, C. and Dykstra, R. (1999) 
Asymptotic behavior of the Grenander estimator at density flat regions.
\newblock \emph{Canad. J. Statist.}, \textbf{27}, 557-566.

\bibitem[{Chen and Samworth(2013)}]{ChenSamworth2013}
Chen, Y. and Samworth, R. J. (2013)
Smoothed log-concave maximum likelihood estimation with applications. 
\newblock \emph{Statist. Sinica.}, \textbf{23}, 1373-1398.

\bibitem[{Chen and Samworth(2014)}]{ChenSamworth2014}
Chen, Y. and Samworth, R. J. (2014)
Generalised additive and index models with shape constraints.
\newblock Available at \texttt{http://arxiv.org/abs/1404.2957}.

\bibitem[{Cs\"org\H o and Horvath (1993)}]{CsorgoHorvath1993}
Cs\"org\H o, M. and Horvath, L. (1993) Weighted approximations in probability and statistics.
\newblock John Wiley \& Sons Inc, Chichester, England. 

\bibitem[{Cule and Samworth(2010)}]{CuleSamworth2010}
Cule, M. and Samworth, R. (2010) Theoretical properties of the log-concave maximum likelihood estimator of a multidimensional density.
\newblock \emph{Electron. J. Statist.}, \textbf{4}, 254--270.

\bibitem[{Doss and Wellner(2016)}]{DossWellner2016}
Doss, C. and Wellner, J. A. (2016)
Global rates of convergence of the MLEs of log-concave and s-concave densities.
\newblock \emph{Ann. Statist.}, to appear. 
\newblock Available at \texttt{http://arxiv.org/abs/1306.1438}.

\bibitem[{Dudley (1999)}]{Dudley1999}
Dudley, R. M. (1999) Uniform central limit theorems. 
\newblock Cambridge University Press, Cambridge, U.K.

\bibitem[{D\"umbgen, Rufibach and Wellner(2007)}]{DRW2007}
D\"umbgen, L., Rufibach, K. and Wellner, J. A. (2007) Marshall's lemma for convex density estimation.
\newblock In \emph{Asymptotics: particles, processes and inverse problems: festschrift for Piet Groeneboom} (E. Cator, G. Jongbloed, C. Kraaikamp, R.
Lopuha\"{a} and J.A. Wellner, eds.), 101--107. Institute of Mathematical Statistics, Ohio.

\bibitem[{D\"umbgen, Samworth and Schuhmacher(2011)}]{DSS2011}
D\"umbgen, L., Samworth, R. and Schuhmacher, D. (2011) Approximation by log-concave distributions with applications to regression. 
\newblock \emph{Ann. Statist.}, \textbf{39}, 702-730.

\bibitem[{Gao and Wellner(2009)}]{GaoWellner2009}
Gao, F. and Wellner, J. A. (2009) On the rate of convergence of the maximum likelihood estimator of a $k$-monotone density
\newblock \emph{Science in China Series A: Mathematics}, \textbf{52}, 1--14.


\bibitem[{Groeneboom(1985)}]{Groeneboom1985}
Groeneboom, P. (1985) Estimating a monotone density.
\newblock In \emph{Proceedings of the Berkeley conference in honor of Jerzy Neyman and Jack Kiefer}, 539--555, Wadsworth, Belmont, CA.

\bibitem[{Groeneboom, Jongbloed and Wellner(2001a)}]{GJW2001a}
Groeneboom, P., Jongbloed, G. and Wellner, J. A. (2001a) A canonical process for estimation of convex functions:
the ``invelope'' of integrated Brownian motion $+t^4$.
\newblock \emph{Ann. Statist.}, \textbf{29}, 1620--1652.

\bibitem[{Groeneboom, Jongbloed and Wellner(2001b)}]{GJW2001b}
Groeneboom, P., Jongbloed, G. and Wellner, J. A. (2001b) Estimation of a convex
  function: characterizations and asymptotic theory.
\newblock \emph{Ann. Statist.}, \textbf{29}, 1653--1698.

\bibitem[{Groeneboom, Jongbloed and Wellner(2008)}]{GJW2008}
Groeneboom, P., Jongbloed, G. and Wellner, J. A. (2008) 
The support reduction algorithm for computing nonparametric function estimates in mixture models. 
\newblock \emph{Scand. J. Statist.}, \textbf{35}, 385-399.

\bibitem[{Groeneboom and Pyke(1983)}]{GroeneboomPyke1983}
Groeneboom, P. and Pyke, R. (1983) Asymptotic normality of statistics based on the convex minorants of empirical distribution functions.
\newblock \emph{Ann. Prob.}, \textbf{11}, 328--345.

\bibitem[{Guntuboyina and Sen(2013)}]{GuntuboyinaSen2013}
Guntuboyina, A. and Sen, B. (2013) Covering numbers for convex functions. 
\newblock \emph{IEEE Trans. Inf. Theory}, \textbf{59}, 1957--1965.

\bibitem[{Guntuboyina and Sen(2015)}]{GuntuboyinaSen2015}
Guntuboyina, A. and Sen, B. (2015)  Global risk bounds and adaptation in univariate convex regression.
\newblock \emph{Prob. Theory Related Fields}, to appear. 
\newblock Available at \texttt{http://arxiv.org/abs/1305.1648}.


\bibitem[{Hanson and Pledger(1976)}]{HansonPledger1976}
Hanson, D.L. and Pledger, G. (1976) Consistency in concave regression.
\newblock \emph{Ann. Statist.}, \textbf{4}, 1038--1050.

\bibitem[{Hildreth(1954)}]{Hildreth1954}
Hildreth, C. (1954) Point estimators of ordinates of concave functions.
\newblock \emph{J. Amer. Statist. Assoc.}, \textbf{49}, 598--619.

\bibitem[{Lachal(1997)}]{Lachal1997}
Lachal, A. (1997) Local asymptotic classes for the successive primitives of Brownian motion
\newblock \emph{Ann. Prob.}, \textbf{25}, 1712--1734.


\bibitem[{Jankowski(2014)}]{Jankowski2014}
Jankowski, H. (2014) Convergence of linear functionals of the Grenander estimator under misspecification
\newblock \emph{Ann. Statist.}, \textbf{42}, 625--653.

\bibitem[{Kim and Pollard(1990)}]{KimPollard1990}
Kim, J. and Pollard, D. (1990) Cube root asymptotics.
\newblock \emph{Ann. Statist.}, \textbf{18}, 191--219.

\bibitem[{Kim and Samworth(2014)}]{KimSamworth2014}
Kim, A. K. H. and Samworth, R. J. (2014) Global rates of convergence in log-concave density estimation.
\newblock Available at \texttt{http://arxiv.org/abs/1404.2298}.

\bibitem[{Mammen(1991)}]{Mammen1991}
Mammen, E. (1991) Nonparametric regression under qualitative smoothness assumptions.
\newblock \emph{Ann. Statist.}, \textbf{19}, 741--759.

\bibitem[{Meyer(2003)}]{Meyer2003}
Meyer, M. C. (2003).  A test for linear versus convex regression function using shape-restricted regression.
\newblock \emph{Biometrika}, \textbf{90}, 223--232.

\bibitem[{Meyer(2013)}]{Meyer2013}
Meyer, M. C. (2013).  Semi-parametric additive constrained regression.
\newblock \emph{J. Nonparametric Statist.}, \textbf{25}, 715--743.

\bibitem[{Meyer and Sen(2013)}]{MeyerSen2013}
Meyer, M. C. and Sen, B. (2013). Testing against a linear regression model using ideas from shape-restricted estimation
\newblock Available at \texttt{http://arxiv.org/abs/1311.6849}.

\bibitem[{Rockafellar(1997)}]{Rockafellar1997}
Rockafellar, R. T. (1997) \emph{Convex Analysis}.
\newblock Princeton University Press, Princeton, NJ.

\bibitem[{Shorack and Wellner(1986)}]{ShorackWellner1986}
Shorack, G. R. and Wellner, J. A. (1986) \emph{Empirical processes with applications to statistics}.
\newblock John Wiley \& Sons Inc, New York, NY.

\end{thebibliography}
\end{document}